# CAUSAL AND ANTICAUSAL OPERATORS ASSOCIATED WITH INPUT- STATE /OUTPUT DESCRIPTIONS OF SWITCHED LINEAR DYNAMIC SYSTEMS WITH POINT LAGS

M. De la Sen, IIDP. **Faculty of Science and Technology**. *University of the Basque Country*.

Campus of Leioa (Bizkaia). Aptdo. 644- Bilbao, SPAIN

**Abstract:** This paper investigates the causality properties of a class of linear time-delay systems under constant delays which possess a finite set of distinct linear time- invariant parameterizations (or configurations) which , together with some switching function, conform a linear time-varying switched dynamic system. Explicit expressions are given to define point-wise the causal and anticausal Toeplitz and Hankel operators from the set of switching time instants generated from the switching function. The case of the auxiliary unforced system defined by the matrix of undelayed dynamics being dichotomic (i.e. it has no eigenvalue on the complex imaginary axis) is considered in detail. Stability conditions as well as dual instability ones are discussed for this case which guarantee that the whole system is either stable, or unstable but no configuration of the switched system has eigenvalues within some vertical strip including the imaginary axis. It is proved that if the system is causal and uniformly controllable and observable then it is globally asymptotically Lyapunov stable independent of the delays provided that a minimum residence time in-between consecutive switches is kept or , if all the set of matrices describing the auxiliary unforced delay – free system commute pair-wise.

## 1. Introduction

The stabilization of dynamic systems is a very important issue since it is the first requirement for most of the applications. Powerful techniques for studying the stability of dynamic systems are Lyapunov stability theory and fixed point theory which can be easily extended from the linear time-invariant case to the time-varying one as well as to functional differential equations , as those arising for instance from the presence of internal delays, and to certain classes of nonlinear systems, [1-2]. Dynamic systems which are of increasing interest are the so-called switched systems which consist of a set of individual parameterizations and a switching function which selects along time which parameterization is active during a subsequent time interval. Switched systems are essentially time-varying by nature even if all the individual parameterizations are time-invariant due to the operation mode of the switching function. The major interest of such systems arises from the fact that some real word existing systems are able to change their parameterizations to better adapt to their environments. Another important interest of some of such systems relies on the fact that changes of parameterizations through time can lead to benefits in certain applications, [3-13]. The natural way of modelling these situations lies in the definition of appropriate switched dynamic systems. For instance, the asymptotic stability of Liénard-type equations with Markovian switching is investigated in [4-5]. Also, time- delay dynamic systems are very important in the real life for appropriate modelling of certain biological and ecological systems and they are present in physical processes implying diffusion, transmission, teleoperation, population dynamics, war and peace models etc. (see, for instance, [1-2], [12-18]). Linear switched dynamic systems are a very particular case of the dynamic system proposed in this manuscript. A switched system can result, for instance, from the use of a multi-model scheme, a multi-controller scheme, a buffer system or a multi-estimation scheme ( see, for instance, [3], [19-24], [17], [27-28]). For instance, a (non exhaustive) list of papers with deal with some of these questions related to switched systems follow:

(1) In [15], the problem of delay-dependent stabilization for singular systems with multiple internal and external incommensurate delays is focused on. Multiple memory-less state- feedback controls are designed so that the resulting closed-loop system is regular independent of delays, impulse- free and asymptotically stable.

(2) In [19], the problem of the N- buffer switched flow networks is discussed based on a theorem on positive topological entropy.

(3) In [20], a multi-model scheme is used for the regulation of the transient regime occurring between stable operation points of a tunnel diode-based triggering circuit.

(4) In [21-22], a parallel multi-estimation scheme is derived to achieve close-loop stabilization in robotic manipulators whose parameters are not perfectly known. The multi-estimation scheme allows the improvement of the transient regime compared to the use of a single estimation scheme while achieving at the same time closed-loop stability.

(5) In [23], a parallel multi-estimation scheme allows the achievement of an order reduction of the system prior to the controller synthesis so that this one is of reduced –order (then less complex) while maintaining closed-loop stability.

(6) In [24], the stabilization of switched dynamic systems is discussed through topologic considerations via graph theory.

(7) The stability of different kinds of switched systems subject to delays has been investigated in [11-13], [17], [27-28].

(8) The stability switch and Hopf bifurcation for a diffusive prey-predator system is discussed in [6] in the presence of delay.

(9) A general theory with discussed examples concerning dynamic switched systems is provided in [3].

The dynamic system under investigation is a linear switched system subject to internal point delays and feedback state- dependent impulsive controls which is based on a finite set of time-varying parametrical configurations and switching function which decides which parameterization is active during a time interval as well as the next switching time instant. Explicit expressions for the state and output trajectories are provided together with the evolution operators and the input –state and input output operators undrer zero initial conditions. The causal and anticausal Toeplitz and causal and anticausal Hankel operators are defined explicitly for the case when all the configurations have auxiliary unforced delay-free systems being dichotomic (i.e. with no eigenvalues on the complex imaginary axis), the controls are square-integrable and the input-output operators are bounded. It is proven that if the anticausal Hankel operator is zero independent of the delays and the system is uniformly controllable and uniformly observable independent of the delays then the system is globally asymptotically Lyapunov´s stable independent of the delays. Those results generalize considerably some parallel background ones for the delay-free and switching-free linear time- invariant case, [28]. The paper is organized as follows. Section 2 discusses the various evolution operators valid to build the state-trajectory solutions in the presence of internal delays and switching functions operating over a set of time-invariant prefixed configurations. Stability and instability are discussed from Gronwall´s lemma,[27] for the case when the auxiliary unforced delay-free system possesses only dichotomic time-invariant configurations. Analytic expressions are given to define such operators as well as the input-state and input-output ones under zero initial conditions. Section 3 discusses the input-state and input-output and operators if the input is square-integrable and the state and output are also square-integrable. Related to those operators proved to be bounded under certain condition, the causal and anticausal state-input and state-output Hankel and the causal and anticausal state-input and state-output Toeplitz operators are defined explicitly. The boundedness of the state-input/ output operators is proven if the controls are square-integrable and the matrices of all the active configurations of the auxiliary –delay free system are dichotomic for the given switching function. The causality and anticausality of the switched system are characterized and some relationships between the properties of causality, stability, controllability and observability are also proven.

### 1.1 Notation

$\mathbf{Z}, \mathbf{R}, \mathbf{C}$ are the sets of integer, real and complex numbers, respectively.

$\mathbf{Z}_+$ and $\mathbf{R}_+$ denote the positive subsets of $\mathbf{Z}$, respectively, and $\mathbf{C}_+$ denotes the subset of $\mathbf{C}$ of complex numbers with positive real part.

$\mathbf{Z}_-$ and $\mathbf{R}_-$ denote the negative subsets of $\mathbf{Z}$, respectively, and $\mathbf{C}_-$ denotes the subset of $\mathbf{C}$ of complex numbers with negative real part.

$\mathbf{Z}_{0+} := \mathbf{Z}_+ \cup \{0\}$, $\mathbf{R}_{0+} := \mathbf{R}_+ \cup \{0\}$, and $\mathbf{C}_{0+} := \mathbf{C}_+ \cup \{0\}$

$\mathbf{Z}_{0-} := \mathbf{Z}_- \cup \{0\}$, $\mathbf{R}_{0-} := \mathbf{R}_- \cup \{0\}$, and $\mathbf{C}_{0-} := \mathbf{C}_- \cup \{0\}$

Given some linear space X (usually $\mathbf{R}$ or $\mathbf{C}$) then $C^{(i)}(\mathbf{R}_{0+}, X)$ denotes the set of functions of class $C^{(i)}$. Also, $BPC^{(i)}(\mathbf{R}_{0+}, X)$ and $PC^{(i)}(\mathbf{R}_{0+}, X)$ denote the set of functions in $C^{(i-1)}(\mathbf{R}_{0+}, X)$ which, furthermore, possess bounded piecewise continuous constant or, respectively, piecewise continuous constant i-th derivative on X.

The set of linear operators from the linear space X to the linear space Y are denoted by $\mathbf{L}(X, Y)$ and the Hilbert space of n norm-square Lebesgue integrable real functions on $\mathbf{R}$ is denoted by $\mathbf{L}_2^n \equiv \mathbf{L}_2^n(\mathbf{R})$ and endowed with the inner product $L^2$–norm $\|f\|_{\mathbf{L}_2^n} := \left(\int_{-\infty}^{\infty} \|f(\tau)\|_2^2 d\tau\right)$, $\forall f \in \mathbf{L}_2^n$, where $\|.\|_2$ is the $\ell_2$–vector (or Euclidean) norm and its corresponding induced matrix norm. $\mathbf{L}_2^n[\alpha, \infty)$ the Hilbert space of n norm-square Lebesgue integrable real functions on $[\alpha, \infty) \subset \mathbf{R}$ for a given $\alpha \in \mathbf{R}$ which is endowed with the norm $\|f\|_{\mathbf{L}_2^n[\alpha, \infty)} := \left(\int_{\alpha}^{\infty} \|f(\tau)\|_2^2 d\tau\right)$, $\forall f \in \mathbf{L}_2^n[\alpha, \infty)$.

$\mathbf{L}_{2+}^n := \{f \in \mathbf{L}_2^n : f(t) = 0, \forall t \in \mathbf{R}_-\}$ and $\mathbf{L}_{2-}^n := \{f \in \mathbf{L}_2^n : f(t) = 0, \forall t \in \mathbf{R}_+\}$ are closed subspaces of $\mathbf{L}_{2+}^n := \{f \in \mathbf{L}_2^n : f(t) = 0, \forall t \in \mathbf{R}_-\}$ of respective supports $\mathbf{R}_{0+}$ and $\mathbf{R}_{0-}$. Then, $\mathbf{L}_2^n = \mathbf{L}_{2+}^n \oplus \mathbf{L}_{2-}^n$.

$I_n$ denotes the n-th identity matrix.

$\lambda_{max}(M)$ and $\lambda_{min}(M)$ stand for the maximum and minimum eigenvalues of a definite square real matrix $M = (m_{ij})$.

$\sigma: \mathbf{R}_{0+} \to \overline{N} := \{1, 2, ..., N\}$ is the switching function which defines the parameterization at time t of a switched dynamic system among N possible time-invariant parameterizations. $\sigma_{\tau,t} (:= \sigma \mid [0,t)): [0,t)(\subset \mathbf{R}_{0+}) \to \overline{N}_{\tau,t} \subset \overline{N}$ is the partial switching function with its domain restricted to $[\tau, t]$. $\sigma_t$ is a notational abbreviation of $\sigma_{0,t}$.

The point constant delays are denoted by $h_i \in [0, h]; \forall i \in \overline{q} \cup \{0\}$ are, in general, incommensurate and $h_0 = 0$.

## 2. The dynamic system subject to time-delays

Consider the following class of switched linear time-varying differential dynamic system subject to q internal incommensurate point delays $0 = h_0 < h_1 < h_2 < ... < h_q = h$:

$$\dot{x}(t) = \sum_{i=0}^{q} A_i(t) x(t - h_i) + B(t) u(t) \; ; \quad y(t) = C(t) x(t) + D(t) u(t) \tag{2.1}$$

where $x(t) \in \mathbf{R}^n$, $u(t) \in \mathbf{R}^m$, and $y(t) \in \mathbf{R}^p$ are the state, input (or control) and output (or measurement) vectors, respectively, and

$A_i(t) \in \mathbf{A}_i := \{A_{ij} \in \mathbf{R}^{n \times n} : j \in \overline{N}\}$; $B(t) \in \mathbf{B} := \{B_j \in \mathbf{R}^{n \times m} : j \in \overline{N}\}$

$C(t) \in \mathbf{C} := \{C_j \in \mathbf{R}^{p \times n} : j \in \overline{N}\}$; $D(t) \in \mathbf{D} := \{D_j \in \mathbf{R}^{p \times m} : j \in \overline{N}\}$

where $i \in \overline{q} \cup \{0\} := \{0, 1, 2, ..., N\}$, fulfilling that $A_i(\tau), B_i(\tau), C_i(\tau)$ and $D_i(\tau)$ are piecewise constant such that they are constant either in $(t-T, t]$ or in $[t, t+T)$, for $\forall t \in \mathbf{R}_{0+}$ and some fixed $T \in \mathbf{R}_+$. The system (2.1) has two auxiliary unforced systems which are useful for stability analysis defined as follows:

*The zero-delay auxiliary unforced switched system (2.1):* $\dot{x}(t) = \left( \sum_{i=0}^{q} A_i(t) \right) x(t) \; ; \; y(t) = C(t) x(t)$

is the particular system arising when all the delays of (2.1) are zero, and

*The delay-free unforced auxiliary switched system:* $\dot{x}(t) = A_0(t) x(t - h_i) \; ; \; y(t) = C(t) x(t)$

is the particular system arising when all the matrices describing delayed dynamics in (2.1) are zero.

A well-known important property is that in the case of one single configuration, i.e. the system does not switch among a set of them, the global stability of the above auxiliary systems lead to necessary conditions for stability independent of the delays, [30]. The physical interpretation is that the dynamic system (2.1) is a switched system under some (piecewise constant) switching function $\sigma: \mathbf{R}_{0+} \to \overline{N}$, which generates a strictly ordered sequence of switching time instants $SI_\sigma := \{t_i : t_{i+1} \geq t_i + T, \forall i \in \{1\} \supset N_0 \subset \mathbf{N}, t_1 \in \mathbf{R}_{0+}\}$, and which might be equivalently rewritten, since $A_i(t) = A_{i\sigma(t)}, \forall i \in \overline{N} \cup \{0\}$, $B(t) = B_{\sigma(t)}$, $C(t) = C_{\sigma(t)}$, $D(t) = D_{\sigma(t)}$ via the switching function $\sigma: \mathbf{R}_{0+} \to \overline{N}$, as:

$$\dot{x}(t) = \sum_{i=0}^{q} A_i(t) x(t-h_i) + B(t) u(t) = \left( \sum_{i=0}^{q} A_i(t) \right) x(t) + \sum_{i=1}^{q} \left( A_i(t) x(t-h_i) - A_i(t) x(t) \right) + B(t) u(t)$$

(2.2)

$$y(t) = C(t) x(t) + D(t) u(t) \tag{2.3}$$

where $x : \mathbf{R}_{0+} \cup [-h, 0) \to X \subset \mathbf{R}^n$ is the state – trajectory solution, which is almost everywhere time-differentiable on $\mathbf{R}_{0+}$ and satisfies (2.2), subject to bounded piecewise continuous initial conditions on $[-h, 0)$, i.e. $x = \varphi \in BPC^{(0)}([-h, 0], \mathbf{R}^n)$. It is assumed that $\sigma(t) = j \in \overline{\mathbf{N}}; \forall t \in \mathbf{R}_- \cup [0, t_1)$, $t_1 \in ST_\sigma$, being the first switching instant generated by the switching function $\sigma : \mathbf{R}_{0+} \to \overline{\mathbf{N}}$; i.e. there is a time –invariant parameterization belonging to the given set on $(-\infty, t_1]$. The above assumption has an obvious real meaning for the general cases where the control is nonzero on $\mathbf{R}_-$. The unique mild solution of the state –trajectory solution, which exists on $\mathbf{R}_{0+}$ according to Picard – Lindeloff theorem for any given $\varphi \in BPC^{(0)}([-h, 0], \mathbf{R}^n)$ and any $u \in BPC^{(0)}(\mathbf{R}, \mathbf{R}^m)$, may be calculated on any time interval $[\alpha, t] \subset \mathbf{R}$ on nonzero measure by first decomposing the interval as a disjoint union of connected components defined by its contained sequence of switching time instants as:

$$[\alpha, t] = [\alpha, t_k] \cup \left( \bigcup_{i \in \hat{N}_t(\alpha)} [t_{k+i}, t_{k+i+1}] \right) \cup [t_{k+1+\overline{N}_t(\alpha)}, t] \tag{2.4}$$

where $\hat{N}_t(\alpha) := \overline{N}_t(\alpha) \cup \{0\}$; $\overline{N}_t(\alpha) := \{ i \in \mathbf{N} : SI_\sigma \ni t_i \leq t \}$, $t_{k+i} \in SI_\sigma$; $\forall i \in \overline{N}_t(\alpha)$ and $t_{k+N_t(\sigma)+1} \in SI_\sigma$. Note that $\sigma(t_k^-) = j(t_k) \neq \sigma(t_k^+) = j(t_{k+1}) \in \overline{\mathbf{N}}; \forall t_k, t_{k+1} \in SI_\sigma$. Then, the state trajectory solution is:

$$x(t) = (\mathbf{\Phi} x_h(\alpha))(t) + (\mathbf{\Gamma} u_\alpha)(t)$$
$$= \Phi(t, \alpha) x(\alpha) + \sum_{i=1}^{q} \int_\alpha^t \Phi(t, \tau) A_i(\tau) x(\tau - h_i) d\tau + \int_\alpha^t \Phi(t, \tau) B(\tau) u(\tau) d\tau \quad (2.5b)$$

where, although the evolution operators between any two time instants $\tau, t > \tau$ depends on the corresponding partial switching function $\sigma_{\tau,t}$, the simpler notation $\Phi(t, \tau)$ is preferred instead for $\Phi_{\sigma_{\tau,t}}(t, \tau)$ for the shake of simplicity, This simplified notation criterion will be used when no confusion is expected together with the former one $M_{\sigma(t)} \to M(t)$ for all the matrices of the individual parameterizations. The output trajectory solution is:

$$y(t) = (C \mathbf{\Phi} x_h(\alpha))(t) + (C \mathbf{\Gamma} u_\alpha + D)(t) \tag{2.6a}$$
$$= C(t) \left[ \Phi(t, \alpha) x(\alpha) + \sum_{i=1}^{q} \int_\alpha^t \Phi(t, \tau) A_i(\tau) x(\tau - h_i) d\tau + \int_\alpha^t \Phi(t, \tau) B(\tau) u(\tau) d\tau \right]$$
$$+ D(t) u(t) \tag{2.6b}$$

$\forall t (\geq \alpha), \alpha \in \mathbf{R}$, subject to initial conditions $\varphi \in BPC^{(0)}([-h, 0], \mathbf{R}^n)$, where:

1) $x_h(\alpha)$ is the strip of state-trajectory solution on $[\alpha-h, h]$ which takes values $\varphi(t)$ if $t = \alpha - h < 0$

2) the evolution operator in $\Phi \in L(\mathbf{R}^n \times \mathbf{R}, \mathbf{R}^n)$ is defined point-wise by

$$(\Phi x_h(\alpha))(t) := \Phi(t,\alpha)x(\alpha) + \sum_{i=1}^{q} \int_{\alpha}^{t} \Phi(t,\tau)x(\tau - h_i)d\tau \;;\; \forall t(\geq \alpha), \alpha \in \mathbf{R} \qquad (2.7)$$

so that $(\Phi x_h(0))(t) := \Phi(t,0)x(0) + \sum_{i=1}^{q} \int_{0}^{t} \Phi(t,\tau)A_i(\tau)x(\tau - h_i)d\tau$ is the unforced response in $[0,t]$, where the matrix function $\Phi \in C^{(0)}(\mathbf{R} \times \mathbf{R}, \mathbf{R}^{n\times n})$ is a fundamental matrix of the dynamic differential system which is everywhere differentiable and has almost everywhere continuous time-derivative on $\mathbf{R}$ with bounded discontinuities on the set $SI_\sigma$ and is defined on the interval $[\alpha,t] \subset \mathbf{R}$ as

$$\Phi(t,\alpha) = e^{A_0\left(t_{k+\overline{N}_t(\alpha)+1}\right)\left(t - t_{k+\overline{N}_t(\alpha)+1}\right)} \left(\prod_{i=1}^{k+\overline{N}_t(\alpha)} \left[e^{A_0(t_{k+i})(t_{k+i+1}-t_{k+i})}\right]\right) e^{A_0(\alpha)(t_k - \alpha)}$$

(2.8)

and the above matrix function products are defined to the left, and

3) the input-state and input-output operators in $\Gamma \in L(\mathbf{R}^m \times \mathbf{R}, \mathbf{R}^n)$ and $\Gamma_o \in L(\mathbf{R}^m \times \mathbf{R}, \mathbf{R}^p)$, respectively $\Gamma_o := C_\sigma \Gamma + D_\sigma$, are defined point-wise by

$$(\Gamma u_{\alpha t})(t) := \int_{\alpha}^{t} \Phi(t,\tau)B(\tau)u(\tau)d\tau = \int_{-\infty}^{t} \Phi(t,\tau)B(\tau)u_{\alpha t}(\tau)d\tau$$

$$(\Gamma_o u_{\alpha t})(t) := \int_{\alpha}^{t} C(t)\Phi(t,\tau)B(\tau)u(\tau)d\tau + D(t)u(t)$$

$$= \int_{-\infty}^{t} C(t)\Phi(t,\tau)B(\tau)u_{\alpha t}(\tau)d\tau + D(t)u(t); \; \forall t(\geq \alpha), \alpha \in \mathbf{R}_{0+} \qquad (2.9)$$

where

$$u_{\alpha t}(\tau) := \begin{cases} u(\tau), \forall \tau \in [\alpha, t] \\ 0, \forall \tau \in \mathbf{R} \setminus [\alpha, t] \end{cases} \qquad (2.10)$$

so that

$$(\Gamma u_{\alpha t})(t) := \sum_{i=1}^{q} \int_{\alpha}^{t} \Phi(t,\tau)u(\tau)d\tau \;;\; (C\Gamma u_{\alpha t} + D)(t) := \sum_{i=1}^{q} \int_{\alpha}^{t} C(t)\Phi(t,\tau)B(\tau)u(\tau)d\tau + D(t)u(t)$$

are, respectively, the unforced state and output responses in $[\alpha, t]$. The state and output trajectory solutions (2.5), or (2.6), under (2.7)-(2.9), subject to the output equation in (2.1) are identically defined by:

$$x(t) = Z(t,\alpha)x(\alpha) + \sum_{i=1}^{q} \int_{\alpha-h_i}^{\alpha} Z(t,\tau)x(\tau)d\tau + \int_{-\infty}^{t} Z(t,\tau)B(\tau)u_{\alpha t}(\tau)d\tau \qquad (2.11a)$$

$$y(t) = C(t)\left[Z(t,\alpha)x(\alpha) + \sum_{i=1}^{q} \int_{\alpha-h_i}^{\alpha} Z(t,\tau)x(\tau)d\tau + \int_{-\infty}^{t} Z(t,\tau)B(\tau)u_{\alpha t}(\tau)d\tau\right] + D(t)u(t)$$

(2.11b)

with initial conditions $x = \varphi \in BPC^{(0)}\left([-h, 0], \mathbf{R}^n\right)$, so that $x(0) = \varphi(0)$, $Z(t, \alpha) \in C^{(0)}\left(\mathbf{R} \times \mathbf{R}, \mathbf{R}^{n \times n}\right)$ is an everywhere differentiable matrix function on $\mathbf{R}_+$, with almost everywhere continuous time-derivative except at time instants in $SI_\sigma$, which satisfies:

$$\dot{Z}(t) = \sum_{i=0}^{q} A(t) Z(t - h_i, 0) \qquad (2.12)$$

on $\mathbf{R}_+$ whose unique solution satisfies $Z(t, \alpha) = 0$; $\forall \alpha (< t), t \in \mathbf{R}$, and is defined by

$$Z(t, \alpha) = \Phi(t, \alpha) \left[ I_n + \sum_{i=1}^{q} \int_{\alpha}^{t} \Phi(t, \tau) A_i(\tau) Z(\tau - h_i, \alpha) d\tau \right]; \forall t (\geq \alpha), \alpha \in \mathbf{R} \qquad (2.13)$$

on any time interval $[\alpha, t] \subset \mathbf{R}_{0+}$. Now, take $\alpha = 0$ and consider that the input u (t) is defined on $\mathbf{R}$. Then, the combination of (2.6) with the substitution of (2.11) in the delayed state and output - trajectory solutions yields:

$$x(t) = \Phi(t, 0) \left( I_n + \sum_{i=1}^{q} \int_0^t \Phi(0, \tau) A_i(\tau) Z(\tau - h_i, 0) d\tau \right) x(0)$$

$$+ \sum_{i=1}^{q} \sum_{j=1}^{q} \int_0^t \int_{-h_j}^{0} \Phi(t, \tau) A_i(\tau) Z(\tau - h_i, \gamma) \varphi(\gamma) d\gamma d\tau$$

$$+ \int_{-\infty}^{t} \Phi(t, \tau) B(\tau) u(\tau) d\tau + \sum_{i=1}^{q} \int_0^t \int_{-\infty}^{\tau - h_i} \Phi(t, \tau) A_i(\tau) Z(\tau - h_i, \gamma) B(\gamma) u(\gamma) d\gamma d\tau$$

(2.14a)

$$= \Phi(t, 0) \left( I_n + \sum_{i=1}^{q} \int_0^t \Phi(0, \tau) A_i(\tau) Z(\tau - h_i, 0) d\tau \right) x(0)$$

$$+ \sum_{i=1}^{q} \sum_{j=1}^{q} \int_0^t \int_{-h_j}^{0} \Phi(t, \tau) A_i(\tau) Z(\tau - h_i, \gamma) \varphi(\gamma) d\gamma d\tau$$

$$+ \int_{-\infty}^{t} \Phi(t, \tau) \left[ B(\tau) + \sum_{i=1}^{q} \int_{-\infty}^{t} \Phi(0, \gamma) A_i(\gamma) Z(\gamma - h_i, \tau) B(\tau) (U(\tau) - U(\gamma - h_i)) d\gamma \right] u(\tau) d\tau$$

(2.14b)

$$y(t) = C_{\sigma(t)} \left[ \Phi(t, 0) \left( I_n + \sum_{i=1}^{q} \int_0^t \Phi(0, \tau) A_i(\tau) Z(\tau - h_i, 0) d\tau \right) x(0) \right.$$

$$+ \sum_{i=1}^{q} \sum_{j=1}^{q} \int_0^t \int_{-h_j}^{0} \Phi(t, \tau) A_i(\tau) Z(\tau - h_i, \gamma) \varphi(\gamma) d\gamma d\tau$$

$$\left. + \int_{-\infty}^{t} \Phi(t, \tau) B(\tau) u(\tau) d\tau + \sum_{i=1}^{q} \int_0^t \int_{-\infty}^{\tau - h_i} \Phi(t, \tau) A_i(\tau) Z(\tau - h_i, \gamma) B(\gamma) u(\gamma) d\gamma d\tau \right] + D(t) u(t)$$

$$= C(t) \left[ \Phi(t, 0) \left( I_n + \sum_{i=1}^{q} \int_0^t \Phi(0, \tau) A_i(\tau) Z(\tau - h_i, 0) d\tau \right) x(0) \right.$$

$$+ \sum_{i=1}^{q} \sum_{j=1}^{q} \int_0^t \int_{-h_j}^{0} \Phi(t, \tau) A_i(\tau) Z(\tau - h_i, \gamma) \varphi(\gamma) d\gamma\, d\tau$$

$$\left. + \int_{-\infty}^{t} \Phi(t, \tau) \left[ B(\tau) + \sum_{i=1}^{q} \int_{-\infty}^{t} \Phi(0, \gamma) A_i(\gamma) Z(\gamma - h_i, \tau) B(\tau) (U(\tau) - U(\gamma - h_i)) d\gamma \right] u(\tau) d\tau \right]$$

$$+ D(t) u(t) \qquad (2.15)$$

where $U(t)$ is the unit step (Heaviside) function. The following result is concerned with sufficient conditions of asymptotic stability and exponential stability of the switched delayed system (2.1)–(2.2), based on Gronwall´s lemma, which will be then useful to define the Hankel and Toeplitz operators.

**Theorem 2.1**. The following properties hold:

**(i)** The unforced dynamic system (2.1)-(2.2) is globally asymptotically stable independent of the sizes of the delays if the switching function $\sigma: \mathbf{R}_{0+} \to \overline{N}$ is such that

$$\exists \lim_{t \to \infty} \prod_{t_i \in ST_\sigma(t)} \left[ K_{0\sigma(t_i)} \left( 1 + \frac{e^{\rho_{0\sigma(t_i)} h} - 1}{\rho_{0\sigma(t_i)}} \left( \sum_{i=1}^{q} \|A_{i\sigma(t_i)}\|_2 \right) \right) e^{-\left( \rho_{0\sigma(t_i)} - K_{0\sigma(t_i)} \sum_{i=1}^{q} \|A_{i\sigma(t_i)}\|_2 \right)(t_{i+1} - t_i)} \right] = 0$$

where $\mathbf{R}_+ \ni K_{0\sigma(t_i)} = K_{0j} (\geq 1) \in \{K_{01}, K_{02}, \dots, K_{0N}\}$ and $\mathbf{R}_+ \ni \rho_{0\sigma(t_i)} \in \{\rho_{01}, \rho_{02}, \dots, \rho_{0N}\}$ if $\sigma(t_i) = j \in \overline{N}$ are real constants such that $\|e^{A_{0i} t}\|_2 \leq K_{0i} e^{-\rho_{0i} t}$; $\forall i \in \overline{N}$ (i.e. all the matrices in the set $\mathbf{A}_0$ are stable) with $ST_\sigma(t) := \{t_i \in ST_\sigma : t_i \leq t\}$ and $t_{s(t)+1} := t$ if $t \notin ST_\sigma$, where $s(t) := \operatorname{card} ST_\sigma(t)$.

**(ii)** The unforced dynamic system (2.1)-(2.2) is globally exponentially stable independent of the sizes of the delays if the switching function $\sigma: \mathbf{R}_{0+} \to \overline{N}$ is such that $A_{0j}$ are all stable matrices satisfying $\rho_{0j} > K_{0\sigma(t_i)} \sum_{i=1}^{q} \|A_{ij}\|_2$ $\forall j \in \overline{N}$, and the residence time at each switching instant satisfies $\max_{t_i \in ST_\sigma} (t_{i+1} - t_i) \geq T$ with its lower-bound $T$ being sufficiently large according to the respective absolute values $|\rho_{0j}|$ of the stability (or convergence) abscissas of $A_{0j}$ (i.e. $-\rho_{0j} < 0$ if all the eigenvalues of $A_{0j}$ are distinct and $-\rho_{0j} + \varepsilon, \varepsilon \to 0^+$, otherwise); $\forall j \in \overline{N}$ and the norms of the matrices $A_{ij}$ $(\forall i \in \overline{q}, j \in \overline{N})$.

**(iii)** The unforced dynamic system (2.1)-(2.2) is globally exponentially stable independent of the sizes of the delays if the switching function $\sigma: \mathbf{R}_{0+} \to \overline{N}$ is such that at least one $A_{0j}$ is a stable matrix satisfying $\rho_{0j} > K_{0\sigma(t_i)} \sum_{i=1}^{q} \|A_{ij}\|_2$, and furthermore, $\max_{t_i, t_{i+1} \in ST_\sigma} (t_{i+1} - t_i : \sigma(t_i) = j)$ is sufficiently

large compared to $\sum_{t_i, t_{i+1} \in ST_\sigma} \max \left( t_{i+1} - t_i : \sigma(t_i) \neq j, \sigma(t_{i+1}) \neq j \right)$, according to the constants $K_{0j}$ $\left( \forall j \in \overline{N} \right)$, the absolute values of the stability abscissas of $A_{0k}$ $\left( \forall k \in \overline{N} \right)$ and norms of $A_{ij}$ $\left( \forall i, j \in \overline{N} \right)$. If there is only a stable matrix $A_{0j}$ in the set $\mathbf{A}_0$. If there is a unique stable matrix $A_{0j}$, for some $j \in \overline{N}$, then the switched system is globally exponentially stable only if the switching function is such that $\sum_{t_k, t_{k+1} \in ST_\sigma} \left( t_{k+1} - t_k : \sigma(t_k) = j \right)$ has infinite measure. If there is a unique stable matrix $A_{0j}$ for some $j \in \overline{N}$ and if the sequence of switching instants $ST_\sigma$ is finite and there is the switching function is such that $\sigma(t_k) = j$ for the last switching instant $t_k$.

**(iv)** If $\mathbf{A}_0 = \mathbf{A}_{0-} \cup \left( \mathbf{A}_{0+} \cup \mathbf{A}_{0\pm} \right)$ where $\mathbf{A}_{0-} \neq \varnothing$, $\mathbf{A}_{0+}$ and $\mathbf{A}_{0\pm}$ are the sets of stable, unstable and critically stable matrices in the set $\mathbf{A}_0$ then the switched system is globally exponentially stable independent of the sizes of the delays if the switching function $\sigma: \mathbf{R}_{0+} \to \overline{N}$ is such that

$$\sum_{t_i, t_{i+1} \in ST_\sigma} \left( t_{i+1} - t_i : \sigma(t_i) = j, A_{0j} \in \mathbf{A}_{0-} \right)$$ is sufficiently large compared to

$$\sum_{t_i, t_{i+1} \in ST_\sigma} \left( t_{i+1} - t_i : \sigma(t_i) = j, \sigma(t_{i+1}) = k, A_{0j}, A_{0k} \notin \mathbf{A}_{0-} \right)$$ according to the constants $K_{0j}$ $\left( \forall j \in \overline{N} \right)$, the absolute values of the stability abscissas of $A_{0k}$ $\left( \forall k \in \overline{N} \right)$ and norms of $A_{ij}$ $\left( \forall i, j \in \overline{N} \right)$.

**Proof**: **(i)** One gets from (2.6) by using Gronwall´s lemma [27]

$$\|x(t)\|_2 \leq \prod_{t_i \in ST_\sigma(t)} \left[ K_{0\sigma(t_i)} \left( 1 + \frac{e^{\rho_{0\sigma(t_i)} h} - 1}{\rho_{0\sigma(t_i)}} \left( \sum_{i=1}^{q} \|A_{i\sigma(t_i)}\|_2 \right) \right) e^{-\left( \rho_{0\sigma(t_i)} - K_{0\sigma(t_i)} \sum_{i=1}^{q} \|A_{i\sigma(t_i)}\|_2 \right) (t_{i+1} - t_i)} \right] \sup_{-h \leq \tau \leq 0} \|\varphi(\tau)\|_2$$

then property (i) follows by simple inspection that it is guaranteed that $\|x(t)\|_2 \to 0$ as $t \to \infty$ since the function of initial condition is bounded on its definition domain.

**(ii)** It follows directly from the above formula since the upper- bounding function of $\|x(t)\|_2$ is of exponential order with decay rate $(-\rho_{0j}) < 0$, provided that $\rho_{0j} > K_{0\sigma(t_i)} \sum_{i=1}^{q} \|A_{ij}\|_2$ $\forall j \in \overline{N}$ provided that the minimum residence time $\max_{t_i \in ST_\sigma} (t_{i+1} - t_i) \geq T$ is sufficiently large. Properties **(iii)** and **(iv)** are direct extensions of Property **(ii)** for the cases when only one delay–free matrix of dynamics is stable or when only a nonempty subset of them are stable matrices, respectively. □

Theorem 2.1 extends known previous ones concerning asymptotic stability of the switched system if all the matrices of the set $\mathbf{A}_0$ are stable and the switching function is subject to a sufficiently large residence time in-between any two consecutive switches. A dual result to Theorem 2.1[(i)-(iii)] is Theorem 2.2 below for instability when all the matrices in the set $\mathbf{A}_0$ are unstable with no stable or critically stable eigenvalues (i.e. all the matrices $A_{0j}; \forall j \in \overline{N}$ are antistable) and the absolute convergence abscissas of $-(A_{0j}); \forall j \in \overline{N}$ are sufficiently large compared to the norms of the matrices of delayed dynamics. Note that although the matrices of delay – free dynamics be antistable, any of the parameterizations of the whole delayed system (2.1)–(2.2) can be antistable since it is well –known that any time-invariant delayed system possessing a principal term in its characteristic polynomial has any unstable value at finite distance and there exists only a finite number of modes within each vertical strip. As a result, the number of unstable eigenvalues is finite and since the system possesses infinitely many eigenvalues, [28], one concludes that the system cannot be antistable.

**Theorem 2.2**. The following properties hold:

**(i)** The unforced dynamic system (2.1)-(2.2) is globally unstable independent of the sizes of the delays if the switching function $\sigma: \mathbf{R}_{0+} \to \overline{N}$ is such that

$$\exists \lim_{t \to \infty} \prod_{t_i \in ST_\sigma(t)} \left[ \left| \tilde{K}_{0\sigma(t_i)} - K_{0\sigma(t_i)} \frac{e^{\rho_{0\sigma(t_i)}h} - 1}{\rho_{0\sigma(t_i)}} \right| \left( \sum_{i=1}^{q} \|A_{i\sigma(t_i)}\|_2 \right) e^{\left( |\tilde{\rho}_{0\sigma(t_i)}| - K_{0\sigma(t_i)} \sum_{i=1}^{q} \|A_{i\sigma(t_i)}\|_2 \right)(t_{i+1} - t_i)} \right] = \infty$$

where $\mathbf{R}_+ \ni \tilde{K}_{0\sigma(t_i)} \leq K_{0\sigma(t_i)} \in \{\tilde{K}_{01}, \tilde{K}_{02}, ..., \tilde{K}_{0N}\}$ and $\mathbf{R}_- \ni \tilde{\rho}_{0\sigma(t_i)} \in \{\tilde{\rho}_{01}, \tilde{\rho}_{02}, ..., \tilde{\rho}_{0N}\}$, with $|\tilde{\rho}_{0j}| \leq |\rho_{0j}|$ (with $\tilde{K}_{0j} \leq K_{0j}$ and $\tilde{\rho}_{0j}$ being located or close to the minimum real part of the eigenvalues of $A_{0j}$ and $\rho_{0j}; \forall j \in \overline{N}$ defined in Theorem 2.1) if $\sigma(t_i) = j \in \overline{N}$ are real constants such that $\|e^{A_{0i}t}\|_2 \geq \tilde{K}_{0i} e^{|\tilde{\rho}_{0i}|t}; \forall i \in \overline{N}$ (i.e. all the matrices in the set $\mathbf{A}_0$ are antistable and then unstable) with $ST_\sigma(t) := \{t_i \in ST_\sigma : t_i \leq t\}$ and $t_{s(t)+1} := t$ if $t \notin ST_\sigma$, where $s(t) := \text{card } ST_\sigma(t)$.

**(ii)** The unforced dynamic system (2.1)-(2.2) is globally exponentially unstable independent of the sizes of the delays if the switching function $\sigma: \mathbf{R}_{0+} \to \overline{N}$ is such that $A_{0j}$ are all unstable matrices satisfying $|\tilde{\rho}_{0j}| > K_{0\sigma(t_i)} \sum_{i=1}^{q} \|A_{ij}\|_2 \ \forall j \in \overline{N}$, and the residence time at each switching instant satisfies $\max_{t_i \in ST_\sigma}(t_{i+1} - t_i) \geq T$ with its lower- bound $T$ being sufficiently large according to the respective absolute values $|\rho_{0j}|$ of the stability abscissas of the stable matrices $(-A_{0j})$ (i.e. $-|\rho_{0j}| < 0$ if all the eigenvalues of $A_{0j}$ are distinct of positive real parts and $-|\tilde{\rho}_{0j}| + \varepsilon \leq -|\rho_{0j}| + \varepsilon, \varepsilon \to 0^+$, otherwise); $\forall j \in \overline{N}$ and norms of $A_{ij} (\forall i, j \in \overline{N})$.

**(iii)** The unforced dynamic system (2.1)-(2.2) is globally exponentially unstable independent of the sizes of the delays if the switching function $\sigma: \mathbf{R}_{0+} \to \overline{N}$ is such that at least one $A_{0j}$ is a stable matrix satisfying $|\rho_{0j}| > K_{0\sigma(t_i)} \sum_{i=1}^{q} \|A_{ij}\|_2$, and furthermore, $\max_{t_i, t_{i+1} \in ST_\sigma}(t_{i+1} - t_i : \sigma(t_i) = j)$ is sufficiently large

compared to $\sum_{t_i, t_{i+1} \in ST_\sigma} \max \left( t_{i+1} - t_i : \sigma(t_i) \neq j, \sigma(t_{i+1}) \neq j \right)$, according to the constants $K_{0j}$ $(\forall j \in \overline{N})$, the absolute values of the stability abscissas of $A_{0k}$ $(\forall k \in \overline{N})$ and norms of $A_{ij}$ $(\forall i, j \in \overline{N})$. If there is only a stable matrix $A_{0j}$ in the set $\mathbf{A}_0$. If there is a unique stable matrix $A_{0j}$, for some $j \in \overline{N}$, then the switched system is globally exponentially stable only if the switching function is such that $\sum_{t_k, t_{k+1} \in ST_\sigma} \left( t_{k+1} - t_k : \sigma(t_k) = j \right)$ has infinite measure. If there is a unique stable matrix $A_{0j}$ for some $j \in \overline{N}$ and if the sequence of switching instants $ST_\sigma$ is finite and there is the switching function is such that $\sigma(t_k) = j$ for the last switching instant $t_k$. □

A combination of Theorems 2.1-2.2 will be used in Section 3 to guarantee the boundedness of the input-state and the input-output operators of the switched system. The following result is direct from the fact that if the system is exponentially stable then its Euclidean norm possesses an upper- bound of exponential order with negative decay rate so that the state and output trajectory solutions are in $\mathbf{L}_2^n$ and $\mathbf{L}_2^p$, respectively. As a result, the input-state $\Gamma$ and input-output $\Gamma_0$ operators are members of $\mathbf{L}\left(\mathbf{L}_2^m, \mathbf{L}_2^n\right)$ and $\mathbf{L}\left(\mathbf{L}_2^m, \mathbf{L}_2^p\right)$, respectively; i.e. linear and then bounded.

**Proposition 2.3**. If any of the properties of Theorem 1[(i)-(iii)] hold for a given switching function $\sigma: \mathbf{R}_{0+} \to \overline{N}$ then the unforced state and output trajectory solutions $(\Phi x_h(\alpha)) \in \mathbf{L}_2^n[\alpha, \infty)$ and $(C\Phi x_h(\alpha)) \in \mathbf{L}_2^p[\alpha, \infty)$; $\forall \alpha \in \mathbf{R}_{0+}$, respectively, so that $\Phi \in \mathbf{L}\left(\mathbf{R}^n \times [-h, 0], \mathbf{L}_2^n[\alpha, \infty)\right)$, and $(C\Phi) \in \mathbf{L}\left(\mathbf{R}^n \times [-h, 0], \mathbf{L}_2^p[\alpha, \infty)\right)$ which are then linear bounded operators since the switched system is either globally asymptotically stable or globally exponentially stable. In particular, $\Phi \in \mathbf{L}\left(\mathbf{R}^n \times [-h, 0], \mathbf{L}_2^n\right)$ and $(C\Phi) \in \mathbf{L}\left(\mathbf{R}^n \times [-h, 0], \mathbf{L}_2^p\right)$.

If, in addition, $u \in \mathbf{L}_2^m[\alpha, \infty)$ for some $\alpha \in \mathbf{R}$ then the respective forced solutions fulfil $(\Gamma u_\alpha) \in \mathbf{L}_2^n[\alpha, \infty)$ and $(\Gamma_o u) \in \mathbf{L}_2^p[\alpha, \infty)$; $\forall \alpha \in \mathbf{R}_{0+}$ which are then bounded operators. Thus, $\Gamma \in \mathbf{L}\left(\mathbf{L}_2^m[\alpha, \infty), \mathbf{L}_2^n[\alpha, \infty)\right)$ and $\Gamma_o \in \mathbf{L}\left(\mathbf{L}_2^m[\alpha, \infty), \mathbf{L}_2^p[\alpha, \infty)\right)$.

If $u \in \mathbf{L}_{2+}^m$ then the respective forced solutions fulfil $(\Gamma_+ u) \in \mathbf{L}_{2+}^n$ and $(\Gamma_{o+} u) \in \mathbf{L}_{2+}^p$, respectively so that $\Gamma_+ \in \mathbf{L}\left(\mathbf{L}_{2+}^m, \mathbf{L}_{2+}^n\right)$ and $\Gamma_{o+} \in \mathbf{L}\left(\mathbf{L}_{2+}^m, \mathbf{L}_{2+}^p\right)$. Equivalently, if $u \in \mathbf{L}_2^m \cap BPC^{(0)}\left(\mathbf{R}_{0+}, \mathbf{R}^m\right)$, i.e. $u(t) = 0, \forall t \in \mathbf{R}_-$, then $\Gamma_+ \in \mathbf{L}\left(\mathbf{L}_{2+}^m, \mathbf{L}_{2+}^n\right)$ and $\Gamma_{o+} \in \mathbf{L}\left(\mathbf{L}_{2+}^m, \mathbf{L}_{2+}^p\right)$. Equivalently, if $u \in \mathbf{L}_2^m \cap BPC^{(0)}\left(\mathbf{R}_{0+}, \mathbf{R}^m\right)$ then $\Gamma \in \mathbf{L}\left(\mathbf{L}_2^m, \mathbf{L}_2^n\right)$ and $\Gamma_{o+} \in \mathbf{L}\left(\mathbf{L}_2^m, \mathbf{L}_2^p\right)$.

If Theorem 2.2 holds and $u \in L_{2-}^m$ then the respective forced solutions fulfil $(\Gamma_- u) \in L_{2-}^n$ and $(\Gamma_{o-} u) \in L_{2-}^p$ so that $\Gamma_- \in L(L_{2-}^m, L_{2-}^n)$ and $\Gamma_{o-} \in L(L_{2-}^m, L_{2-}^p)$. Equivalently, if $u \in L_2^m \cap BPC^{(0)}(\mathbf{R}_{0-}, \mathbf{R}^m)$ then $\Gamma \in L(L_2^m, L_2^n)$ and $\Gamma_{o+} \in L(L_2^m, L_2^p)$.

**Proof**: The first part concerning the unforced solution follows directly from Theorem 2.1 (i)–(iii). The respective linear operators are bounded. The second part follows by taking into account the above properties in Theorem 2.1 and the square-integrability of u on its appropriate definition domains. □

If the system is globally asymptotically stable then it is possible to restrict the domain and image $L_2^m$ and $L_2^n$ of $\Gamma_+$ to $L_{2+}^m$ and $L_{2+}^n$, respectively, for vector functions $u \in L_{2+}^m$ such that $(\Gamma_+ u) \in L_{2+}^n$, since their support is $\mathbf{R}_{0+}$ and $(\Gamma u)(t) = 0, \forall t \in \mathbf{R}_{0-}$, and then to define a restricted operator $\Gamma_+ \in L(L_{2+}^m, L_{2+}^n)$. In the same way, it is possible to define a restricted operator $(C\Gamma + D)_+ \in L(L_{2+}^m, L_{2+}^p)$. Similarly, it is possible to define $\Gamma_- \in L(L_{2-}^m, L_{2-}^n)$ and $(C\Gamma + D)_- \in L(L_{2-}^m, L_{2-}^p)$ of usefulness for vector functions $u \in L_{2-}^m$ if the system is unstable. □

## 3. Input-state and input-to-output operators of the switched system and Hankel and Toeplitz operators $\mathbf{R}_y \equiv \mathbf{R}_\pm$

This section investigates the input-state and input-output operators $\Gamma_{[0,t]}: \mathbf{R}^m \times [0,t] \to \mathbf{R}^n \times [0,t]$ and $\Gamma_{o[0,t]}: \mathbf{R}^m \times [0,t] \to \mathbf{R}^p \times [0,t]$ of the switched system (2.1) and explicit expressions defining them are given. Then, if the input is a square-integrable real m-vector on $\mathbf{R}$, further conditions for $\Gamma \in L(L_2^m, L_2^n)$ and $\Gamma_o \in L(L_2^m, L_2^p)$ are investigated and weaker ones are also given for $\Gamma_{\mathbf{R}_y} \in L(L_2^m[\mathbf{R}_y], L_2^n[\mathbf{R}_y])$ o $\Gamma_{o\mathbf{R}_y} \in L(L_2^m[\mathbf{R}_y], L_2^p[\mathbf{R}_y])$ with $\mathbf{R}_y \subset \mathbf{R}$ being a bounded real interval, in particular for $\mathbf{R}_y \equiv \mathbf{R}_\pm$. Finally, The Hankel and Toeplitz causal and anticausal operators are investigated concerning the cases $\mathbf{R}_y \equiv \mathbf{R}_\pm$. Two different sets of assumptions, the first one less restrictive, are now given to be used when deriving some of the results of this section:

**Assumptions 3.1**. $u \in L_2^m \cap BPC^{(0)}(\mathbf{R}, \mathbf{R}^m)$, the matrices $A_{0j}$ are dichotomic (i.e. they have no eigenvalues on the imaginary axis) while they have stable and antistable diagonal blocks $A_{0j}^-$ and $A_{0j}^+$ of the same respective orders $n^-$ and $n^+$; $\forall j \in \overline{N}$ which satisfy $\min_{j \in \overline{N}}(|\rho_{0j}|, |\tilde{\rho}_{0j}|) \geq \varepsilon \in \mathbf{R}_+$. Furthermore, the norms of all the matrices of delayed dynamics are less than $\varepsilon$ so that Theorem 2.1 (respectively, Theorem 2.2) holds if all the matrices in the set $A_0$ are stable (respectively, antistable).

□

**Assumptions 3.2**. Assumptions 3.1 hold and, furthermore, the matrices $A_{0j}$ are simultaneously block diagonalizable through the same transformation matrix; $\forall j \in \overline{N}$. □

Note that if Assumptions 3.1 hold there no configuration of the switched system has eigenvalues within the open vertical strip $(-\varepsilon, \varepsilon) \times \mathbf{R}$ of the complex plane from Theorems 2.1-2.2. Furthermore, there exist non unique coordinate transformations $T_j$; $\forall j \in \overline{N}$ such that:

$$\tilde{A}_{0j} = T_j A_{0j} T_j^{-1} = \text{Block Diag}\left[\tilde{A}_{0j}^{-}, \tilde{A}_{0j}^{+}\right] \tag{3.1}$$

where $\tilde{A}_{0j-}$ is stable (i.e. all its eigenvalues are in $\mathbf{C}_{0-}$) and of order $n^-$ and $\tilde{A}_{0j+}$ is antistable (i.e. all its eigenvalues are in $\mathbf{C}_{0+}$) and of order $n^+ = n - n^-$; $\forall j \in \overline{N}$. Note also that if Assumptions 3.2 hold then $T_j = T$, $\forall j \in \overline{N}$. After a linear change of variables $\tilde{x}(t) = T(t_k) x(t)$; $\forall t \in [t_k, t_{k+1})$ with $t_k, t_{k+1} \in ST_\sigma$, such that $\sigma(t_k) = j$ and $T(t_k) = T_j$, for some $j \in \overline{N}$, the system (2.1) may be described as follows:

$$\dot{\tilde{x}}(t) = \sum_{i=0}^{q} \tilde{A}_i(t) \tilde{x}(t - h_i) + \tilde{B}(t) u(t) \quad ; \quad y(t) = \tilde{C}(t) \tilde{x}(t) + \tilde{D}(t) u(t) \tag{3.2}$$

$\forall t \in [t_k, t_{k+1})$, where

$$\tilde{A}_i(t) = \tilde{A}_{ik} \in \tilde{\mathbf{A}}_i := \{\tilde{A}_{ij} \in \mathbf{R}^{n \times n} : j \in \overline{N}\}; \quad \tilde{B}(t) = \tilde{B}_k \in \tilde{\mathbf{B}} := \{\tilde{B}_j \in \mathbf{R}^{n \times m} : j \in \overline{N}\}$$

$$\tilde{C}(t) = \tilde{C}_k \in \tilde{\mathbf{C}} := \{\tilde{C}_j \in \mathbf{R}^{p \times n} : j \in \overline{N}\}; \quad \tilde{D}(t) = \tilde{D}_k \in \tilde{\mathbf{D}} := \{\tilde{D}_j \in \mathbf{R}^{p \times m} : j \in \overline{N}\} \tag{3.3}$$

for some $k \in \overline{N}$ such that $\sigma(t) = k$ subject to (3.1) and (3.4) below:

$$\tilde{A}_i(t) = \tilde{A}_{ij} = T_j A_{ij} T_j^{-1} = \begin{bmatrix} \tilde{A}_{ij}^{--} & \tilde{A}_{ij}^{-+} \\ \tilde{A}_{ij}^{+-} & \tilde{A}_{ij}^{++} \end{bmatrix} ; \quad \forall i \in \overline{q} \cup \{0\} \text{ and some } \forall j \in \overline{N} \tag{3.4}$$

$$\tilde{B}(t) = \tilde{B}_j = T_j B_j = \begin{bmatrix} \tilde{B}_j^{-} \\ \tilde{B}_j^{+} \end{bmatrix} = \tilde{C}(t) = \tilde{C}_j = T_j C_j T_j^{-1} = \text{Block Diag}\left[\tilde{C}_j^{-}, \tilde{C}_j^{+}\right], \tilde{D}(t) = \tilde{D}_j = D_j \tag{3.5}$$

; $\forall i \in \overline{q} \cup \{0\}$ and some $\forall j \in \overline{N}$, where :

$\tilde{A}_{0j}^{-} \in \mathbf{R}^{n^- \times n^-}$, $\tilde{A}_{0j}^{-} \in \mathbf{R}^{n^+ \times n^+}$, $\tilde{A}_{ij}^{--} \in \mathbf{R}^{n^- \times n^-}$, $\tilde{A}_{ij}^{-+} \in \mathbf{R}^{n^- \times n^+}$, $\tilde{A}_{ij}^{+-} \in \mathbf{R}^{n^+ \times n^-}$ and $\tilde{A}_{ij}^{++} \in \mathbf{R}^{n^+ \times n^+}$, $\tilde{B}_j^{-} \in \mathbf{R}^{n^- \times m}$, $\tilde{B}_j^{+} \in \mathbf{R}^{n^+ \times m}$, $\tilde{C}_j^{-} \in \mathbf{R}^{p \times n^-}$ and $\tilde{C}_j^{+} \in \mathbf{R}^{p \times n^+}$, $\forall i \in \overline{q} \cup \{0\}, \forall j \in \overline{N}$.

The subspaces $\chi^{-}(A_{0j}) = \text{Im } T_j^{-1} \begin{bmatrix} I_{n^-} \\ 0 \end{bmatrix}$ and $\chi^{+}(A_{0j}) = \text{Im } T_j^{-1} \begin{bmatrix} 0 \\ I_{n^+} \end{bmatrix}$ are independent of $T_j$ and are called, respectively, the stable and antistable subspaces of $A_{0j}$; $\forall j \in \overline{N}$ which are complementary, i.e.

$\mathbf{R}^n = \chi^-(A_{0j}) \oplus \chi^+(A_{0j})$; $\forall j \in \overline{N}$ so that $\mathbf{R}^n = \chi^-(A_0(t)) \oplus \chi^+(A_0(t))$; $\forall t \in \mathbf{R}_0^+$. The projections on those subspaces are given by the respective formulas:

$$\Pi_j^- = T_j^{-1}\begin{bmatrix} I_{n_-} & 0 \\ 0 & 0 \end{bmatrix} T_j \text{ and } \Pi_j^+ = T_j^{-1}\begin{bmatrix} 0 & 0 \\ 0 & I_{n_+} \end{bmatrix} T_j; \forall j \in \overline{N} \quad (3.6a)$$

$$\Pi^-(t) = T^{-1}(t)\begin{bmatrix} I_{n_-} & 0 \\ 0 & 0 \end{bmatrix} T(t) = T_j^{-1}\begin{bmatrix} I_{n_-} & 0 \\ 0 & 0 \end{bmatrix} T_j; \Pi^+(t) = T^{-1}(t)\begin{bmatrix} I_{n_-} & 0 \\ 0 & 0 \end{bmatrix} T(t) = T_j^{-1}\begin{bmatrix} I_{n_-} & 0 \\ 0 & 0 \end{bmatrix} T_j$$
(3.6b)

and $\Pi_{j_k}^\pm = \Pi^\pm(t) = \Pi^\pm(t_k)$ for some $j_k \in \overline{N}$ such that $\sigma(t) = j_k$; $\forall t \in [t_k, t_{k+1})$ for each $t_k, t_{k+1} \in ST(\sigma)$. Thus, from (2.11)–(2.13) and $\alpha = 0$, one gets directly:

$$\dot{\tilde{x}}(t) = T(t)x(t) = \tilde{Z}(t,0)T(0)x(0) + \sum_{i=1}^{q}\int_{-h_i}^{0} \tilde{Z}(t,\tau)T(\tau)x(\tau)d\tau + \int_{-\infty}^{t} \tilde{Z}(t,\tau)\tilde{B}(\tau)u(\tau)d\tau \quad (3.7)$$

with initial conditions $\tilde{x} = \tilde{\varphi} = T(0)\varphi \in BPC^{(0)}([-h, 0], \mathbf{R}^n)$, so that $\tilde{x}(0) = \tilde{\varphi}(0)$, $\tilde{Z}(t,0) \in C^{(0)}(\mathbf{R}_{0+}, \mathbf{R}^{n \times n})$ is an everywhere differentiable matrix function on $\mathbf{R}_{0+}$, with almost everywhere continuous time-derivative except at time instants in $SI_\sigma$, which satisfies:

$$\dot{\tilde{Z}}(t) = T(t_k)\dot{Z}(t) = \sum_{i=0}^{q} \tilde{A}_i(t_k)T(t-h_i)Z(t-h_i, 0); \forall t \in [t_k, t_{k+1}), \forall t_k \in ST(\sigma) \quad (3.8)$$

on $\mathbf{R}_+$, since $\sigma(t) = \sigma(t_k)$; $\forall t \in [t_k, t_{k+1})$ whose unique solution satisfies $Z(t, 0) = 0$, $\forall t \in \mathbf{R}_-$, and is defined by

$$\tilde{Z}(t,0) = T(t)Z(t,0) = \tilde{\Phi}(t,0)\left[I_n + \sum_{i=1}^{q}\int_0^t \tilde{\Phi}(t,\tau)\tilde{A}_i(\tau)T(\tau-h_i)Z(\tau-h_i, 0)d\tau\right]; \forall t \in \mathbf{R}_{0+}$$
(3.9)

Then,

$$e^{\tilde{A}_0(t)t} = T(t)e^{A_0(t)t}T^{-1}(t) = \text{Block Diag}\left[e^{\tilde{A}_0^-(t)t}, e^{\tilde{A}_0^+(t)t}\right] \quad (3.10a)$$

so that

$$\tilde{\Phi}(t, t_k) = e^{\tilde{A}_0(t_k)(t-t_k)} = T(t_k)\Phi(t, t_k)T^{-1}(t_k)$$

$$= \text{Block Diag}\left[\tilde{\Phi}^-(t, t_k), \tilde{\Phi}^+(t, t_k)\right] = \text{Block Diag}\left[e^{\tilde{A}_0^-(t_k)(t-t_k)}, e^{\tilde{A}_0^+(t_k)(t-t_k)}\right] \quad (3.10b)$$

$$\tilde{Z}(t, t_k) = T(t)Z(t, t_k)T^{-1}(t) = T_{\sigma(t)}Z(t, t_k)T_{\sigma(t)}^{-1} \quad (3.10c)$$

$\forall t \in [t_k, t_{k+1})$, with $\sigma(t) = \sigma(t_k) = j_k \in \overline{N}$, $\forall t \in [t_k, t_{k+1})$, $\forall t_k, t_{k+1} \in ST_\sigma$ since

$$e^{\tilde{A}_{0j}t} = T_j e^{A_{0j}t} T_j^{-1} = \text{Block Diag}\left[e^{\tilde{A}_{0j}^- t}, e^{\tilde{A}_{0j}^+ t}\right] \quad (3.11a)$$

$$\tilde{\Phi}(t, t_k) = T_j \Phi(t, t_k) T_j^{-1}; \tilde{Z}(t, t_k) = T_j Z(t, t_k) T_j^{-1} \quad (3.11b)$$

$\forall t \in [t_k, t_{k+1})$, $\forall t_k \in ST(\sigma)$ provided that $\sigma(t) = \sigma(t_k) = j = j_k \in \overline{N}$ and the transformations also apply on the evolution operators when performing the change of variables.

The input-state and input-output operators obtained in (2.11), (2.14) and (2.15), by taking into account (3.5), are now expanded in the subsequent result for a switching function $\sigma: \mathbf{R}_{0+} \to \overline{N}$. Note that the input-state operator depends on the state variable transformations while the input-output operator does not depend on the state variables, i.e. it does not depend on the matrices $T_{\sigma(t)}$.

**Lemma 3.3**. The input-state and input-output operators have the following point-wise expressions:

$$(\Gamma u)(t) = \int_{-\infty}^{t} Z(t,\tau) B(\tau) u(\tau) d\tau$$

$$= \int_{-\infty}^{t} \Phi(t,\tau) \left[ B(\tau) + \sum_{i=1}^{q} \int_{-\infty}^{t} \Phi(0,\gamma) A_i(\gamma) Z(\gamma - h_i, \tau) B(\tau) (U(\tau) - U(\gamma - h_i)) d\gamma \right] u(\tau) d\tau$$

(3.11a)

$$= \int_{-\infty}^{t} Z(t,\tau) \Pi^-(\tau) B(\tau) u(\tau) d\tau - \int_{t}^{\infty} Z(t,\tau) \Pi^+(\tau) B(\tau) u(\tau) d\tau \quad (3.11b)$$

$$= \int_{-\infty}^{t} \Phi(t,\tau) \left[ \Pi^-(\tau) B(\tau) + \sum_{i=1}^{q} \int_{-\infty}^{t} \Phi(0,\gamma) A_i(\gamma) Z(\gamma - h_i, \tau) \Pi^-(\tau) B(\tau) (U(\tau) - U(\gamma - h_i)) d\gamma \right] u(\tau) d\tau$$

$$- \int_{t}^{\infty} \Phi(t,\tau) \left[ \Pi^+(\tau) B(\tau) + \sum_{i=1}^{q} \int_{-\infty}^{t} \Phi(0,\gamma) A_i(\gamma) Z(\gamma - h_i, \tau) \Pi^+(\tau) B(\tau) (U(\tau) - U(\gamma - h_i)) d\gamma \right] u(\tau) d\tau$$

(3.11c)

$$(\Gamma_o u)(t) = \int_{-\infty}^{t} C(t) Z(t,\tau) B(\tau) u(\tau) d\tau + D(t) u(t)$$

$$= \int_{-\infty}^{t} C(t) \Phi(t,\tau) \left[ B(\tau) + \sum_{i=1}^{q} \int_{-\infty}^{t} \Phi(0,\gamma) A_i(\gamma) Z(\gamma - h_i, \tau) B(\tau) (U(\tau) - U(\gamma - h_i)) d\gamma \right] u(\tau) d\tau + D(t) u(t)$$

(3.12a)

$$= \int_{-\infty}^{t} C(t) Z(t,\tau) \Pi^-(\tau) B(\tau) u(\tau) d\tau - \int_{t}^{\infty} C(t) Z(t,\tau) \Pi^+(\tau) B(\tau) u(\tau) d\tau + D(t) u(t) \quad (3.12b)$$

$$= \int_{-\infty}^{t} C(t) \Phi(t,\tau) \left[ \Pi^-(\tau) B(\tau) + \sum_{i=1}^{q} \int_{-\infty}^{t} \Phi(0,\gamma) A_i(\gamma) Z(\gamma - h_i, \tau) \Pi^-(\tau) B_\sigma(\tau) (U(\tau) - U(\gamma - h_i)) d\gamma \right] u(\tau) d\tau$$

$$- \int_{t}^{\infty} C(t) \Phi(t,\tau) \left[ \Pi^+(\tau) B(\tau) + \sum_{i=1}^{q} \int_{-\infty}^{t} \Phi(0,\gamma) A_i(\gamma) Z(\gamma - h_i, \tau) \Pi^+(\tau) B(\tau) (U(\tau) - U(\gamma - h_i)) d\gamma \right] u(\tau) d\tau$$

$$+ D(t) u(t) \quad (3.12c)$$

**Proof**: It follows directly since the forced solutions of (2.14) –(2.15) may be recalculated by direct manipulation of the integrals as follows:

$$\int_{-\infty}^{t} \Phi(t,\tau)B(\tau)u(\tau)d\tau + \sum_{i=1}^{q} \int_{-\infty}^{\tau-h_i} \int_{0}^{t} \Phi(t,\tau) A_i(\tau) Z(\tau-h_i,\gamma)B(\gamma)u(\gamma)d\gamma d\tau$$

$$= \int_{-\infty}^{t} \Phi(t,\tau)B(\tau)u(\tau)d\tau$$

$$+ \sum_{i=1}^{q} \int_{-\infty}^{t} \int_{0}^{t} \Phi(t,\tau) A_i(\tau) Z(\tau-h_i,\gamma)B(\gamma)(U(\gamma)-U(\tau-h_i))u(\gamma)d\tau d\gamma$$

$$= \int_{-\infty}^{t} \Phi(t,\tau)B_{\sigma(\tau)}(\tau)u(\tau)d\tau$$

$$+ \int_{-\infty}^{t} \sum_{i=1}^{q} \int_{0}^{t} \Phi(t,\tau) A_i(\tau) Z(\tau-h_i,\gamma)B(\tau)(U(\gamma)-U(\tau-h_i))u(\gamma)d\tau d\gamma$$

$$= \int_{-\infty}^{t} \Phi(t,\tau)B(\tau)u(\tau)d\tau$$

$$+ \int_{-\infty}^{t} \sum_{i=1}^{q} \int_{0}^{t} \Phi(t,\gamma) A_i(\tau) Z(\gamma-h_i,\tau)B(\tau)(U(\tau)-U(\gamma-h_i))u(\tau) d\gamma d\tau \qquad \square$$

Now (3.11c)-(3.12c) are further expanded by using the transformation of state variables and the contribution of each inter- switching time intervals. The subsequent auxiliary useful notation convention is used to write the mathematical expressions in a very comprehensive way. It is taken into account that there are no switching instants at negative time, that the current time t may be or not to be a switching instant and that the transformation of variables are given by a non-singular matrix $T_{\sigma(t)}$ which takes a finite number of N values and which is constant within the semi-open time interval in-between any two consecutive switching instants:

$$t_0 = -\infty, \; t_i \in \mathbf{R}_{0+} : \begin{cases} t_i \in ST_\sigma, \forall i \in \overline{k(t)} \text{ if } t_{k(t)} = t \in SI_\sigma \\ t_i \in ST_\sigma, \forall i \in \overline{k(t)-1} \text{ if } t_{k(t)} = t \notin SI_\sigma \end{cases} ; \; \forall i \in \overline{k(t)}; \; t_\infty = \infty \qquad (3.13)$$

where $k: \mathbf{R}_{0+} \to \overline{N}$ is a discrete valued function which takes only a finite number of positive integers according to the switching function used.

**Lemma 3.4**. The input-state and input-output operators have the following expressions:

$$(\mathbf{\Gamma} u)(t) =$$

$$\sum_{j=1}^{k(t)} \int_{t_{j-1}}^{t_j} \tilde{\Phi}^-(t,\tau)\left[\tilde{B}^-(t_{j-1}) + \sum_{i=1}^{q} \sum_{\ell=1}^{k(t)} \int_{t_{\ell-1}}^{t_\ell} \left(\tilde{\Phi}(0,\gamma)\tilde{A}(t_{\ell-1}) \tilde{Z}(\gamma-h_i,\tau)\tilde{B}^-(t_{j-1})\right)^- (U(\tau)-U(\gamma-h_i))d\gamma\right]u(\tau)d\tau$$

$$- \sum_{j=k(t)}^{\infty} \int_{t_j}^{t_{j+1}} \tilde{\Phi}^+(t,\tau)\left[\tilde{B}^+(t_{j-1}) + \sum_{i=1}^{q} \sum_{\ell=1}^{k(t)} \int_{t_{\ell-1}}^{t_\ell} \left(\tilde{\Phi}(0,\gamma)\tilde{A}_i(t_{\ell-1}) \tilde{Z}(\gamma-h_i,\tau)\tilde{B}(t_{j-1})\right)^+ (U(\tau)-U(\gamma-h_i))d\gamma\right]u(\tau)d\tau$$

(3.14a)

$$(\Gamma_o u)(t) =$$

$$\sum_{j=1}^{k(t)} \int_{t_{j-1}}^{t_j} \tilde{C}^-(t)\tilde{\Phi}^-(t,\tau)\left[\tilde{B}^-(t_{j-1}) + \sum_{i=1}^{q}\sum_{\ell=1}^{k(t)}\int_{t_{\ell-1}}^{t_\ell}\left(\tilde{\Phi}(0,\gamma)\tilde{A}_i(t_{\ell-1})\tilde{Z}(\gamma-h_i,\tau)\tilde{B}(t_{j-1})\right)^-(U(\tau)-U(\gamma-h_i))d\gamma\right]u(\tau)d\tau$$

$$-\sum_{j=k(t)}^{\infty} \int_{t_j}^{t_{j+1}} \tilde{C}^+(t)\tilde{\Phi}^+(t,\tau)\left[\tilde{B}^-(t_{j-1}) + \sum_{i=1}^{q}\sum_{\ell=1}^{k(t)}\int_{t_{\ell-1}}^{t_\ell}\left(\tilde{\Phi}(0,\gamma)\tilde{A}_i(t_{\ell-1})\tilde{Z}(\gamma-h_i,\tau)\tilde{B}(t_{j-1})\right)^+(U(\tau)-U(\gamma-h_i))d\gamma\right]u(\tau)d\tau$$

$$+ D(t)u(t) \qquad (3.14b)$$

**Proof**: It follows directly from Lemma 3.3 by using (3.5), (3.6) and (3.10b), since (3.9)-(3.10) hold, where:

$$\tilde{\Phi}(0,\gamma)\tilde{A}_i(\gamma)\tilde{Z}(\gamma-h_i,\tau)\tilde{B}(\tau)$$

$$= T(t)\Phi(0,\gamma)A_i(\gamma)Z(\gamma-h_i,\tau)B(\tau) = \begin{bmatrix}\left(\tilde{\Phi}(0,\gamma)\tilde{A}_i(\gamma)\tilde{Z}(\gamma-h_i,\tau)\tilde{B}(\tau)\right)^- \\ \left(\tilde{\Phi}(0,\gamma)\tilde{A}_i(\gamma)\tilde{Z}(\gamma-h_i,\tau)\tilde{B}(\tau)\right)^+\end{bmatrix} \qquad (3.15)$$

$\forall \gamma \in [t_j, t_{j+1})$, $\sigma(\gamma) = \sigma(t_j)$, $t_0 = -\infty$, $t_j \in ST_\sigma$ for $j \in \overline{N}_0 \subset N$. □

Lemmas 3.3-3.4 will be then used for the explicit definition of the Hankel and Toeplitz operators of the input-state and input-output operators. The following result is useful as an auxiliary one for a subsequent specification of Lemmas 3.3-3.4 either for the general case or for the cases when either Assumptions 3.1 or Assumptions 3.2 hold.

**Lemma 3.5.** The following properties hold:

**(i)** $\exists T(t,\tau): \mathbf{R}_0^2 \to \mathbf{R}^{n \times n}$ dependent on the switching instants $ST_\sigma$ on $[\tau,t) \subset \mathbf{R}_{0+}$ (i.e. $T(t,\tau) = T_{\sigma_{\tau,t}}$ depends on the partial switching function $\sigma_{\tau,t}: [\tau,t) \to \overline{N}_{\tau,t} \subset \overline{N}$) such that:

$$\tilde{\Phi}(t,\tau) = \text{Block Diag}\left[\tilde{\Phi}^-(t,\tau), \tilde{\Phi}^+(t,\tau)\right] = T(t,\tau)\Phi(t,\tau)T^{-1}(t,\tau) \qquad (3.16)$$

which is non-singular for any finite arguments irrespective of Assumptions 3.1, where $\sigma_{\tau,t}(:= \sigma \mid [\tau,t)): [\tau,t)(\subset \mathbf{R}_{0+}) \to \overline{N}_{\tau,t} \subset \overline{N}$ is the partial switching function with its domain restricted to $[\tau,t)$. $\tilde{\Phi}^-(t,\tau)$ and $\tilde{\Phi}^+(t,\tau)$ are, in general, of time interval- dependent sizes $n^-(\tau,t), n^+(\tau,t)$, respectively.

**(ii)** If Assumptions 3.1 hold then

$$\tilde{\Phi}(\tau,t_i) = \text{Block Diag}\left[\tilde{\Phi}^-(\tau,t_i), \tilde{\Phi}^+(\tau,t_i)\right] = T(t_i)\Phi(\tau,t_i)T^{-1}(t_i) \qquad (3.17)$$

$$= \text{Block Diag}\left[e^{\tilde{A}_0^-(t_i)(\tau-t_i)}, e^{\tilde{A}_0^+(t_i)(\tau-t_i)}\right], \forall \tau \in [t_i, t_{i+1}); \forall t_i, t_{i+1} \in ST_\sigma \qquad (3.18)$$

with the first and second square matrix function blocks being convergent and divergent, respectively, and of associate time invariant sizes $n^-, n^+$.

**(iii)** If Assumptions 3.2 hold then **(ii)** holds with constant $T_{\sigma(t)} = T$; $\forall t \in \mathbf{R}_{0+}$.

**(iv)** If Assumptions 3.2 hold and all the matrices in the set $\mathbf{A}_0$ defining the switched system by the partial switching function up to time t defined as $\sigma_t \equiv \sigma_{0,t}(:= \sigma|[0,t)):[0,t) \to \overline{N}_t \subset \overline{N}$ commute, so that $T_{\sigma(t)} = T$, then

$$\tilde{\Phi}(t,\tau) = \text{Block Diag}\left[\tilde{\Phi}^-(t,\tau), \tilde{\Phi}^+(t,\tau)\right] = T\Phi(t,\tau)T^{-1} \tag{3.19}$$

$$= \text{Block Diag}\left[ e^{\sum_{i=1}^{k(t)} \tilde{A}_0^-(t_{i-1})(t_i - t_{i-1})}, \; e^{\sum_{i=1}^{k(t)} \tilde{A}_0^+(t_{i-1})(t_i - t_{i-1})} \right] \tag{3.20}$$

, $\forall t, \tau \in \mathbf{R}_{0+}$; $\forall t_i, t_{i+1} \in ST_{\sigma_t}$. Furthermore,

$$\tilde{\Phi}^{\pm}(t,\tau)\left(\tilde{\Phi}(0,\gamma)\tilde{A}_i(\gamma)\tilde{Z}(\gamma-h_i,\tau)\tilde{B}(\tau)\right)^{\pm} = \tilde{\Phi}^{\pm}(t,\gamma)\left(\tilde{A}_i(\gamma)\tilde{Z}(\gamma-h_i,\tau)\tilde{B}(\tau)\right)^{\pm}$$

in (3.13)-(3.14) subject to (3.20).

**(v)** If both assumptions of Property **(vi)** hold and all the matrices in the set $\mathbf{A}_i$ ($\forall i \in \overline{q}$) defining the switched system by the partial switching function up to time t have a block diagonal structure with two block matrices of common sizes $n^-$ and $n^+$ then $Z(t,\tau)$ is block diagonalizable with two non-zero square matrix blocks of time invariant sizes $n^-$ and $n^+$; $\forall \tau \in \mathbf{R}\setminus(t,\infty), \forall t \in \mathbf{R}$. Furthermore,

$$\left(\tilde{\Phi}(t,\tau)\tilde{A}_i(\tau)\tilde{Z}(\gamma-h_i,\tau)\tilde{B}(\tau)\right)^{\pm} = \tilde{\Phi}^{\pm}(t,\tau)\tilde{A}_i^{\pm}(\tau)\tilde{Z}^{\pm}(\gamma-h_i,\tau)\tilde{B}^{\pm}(\tau); \; \forall t, \tau \in \mathbf{R}$$

in (3.13)-(3.14).

**Proof**: **(i)** It follows directly from the fact that any real matrix has a Jordan diagonal form.

**(ii)-(iii)** They follow directly from the fact that the matrix function $\Phi(\tau, t_i)$ is an exponential matrix function of $A_{0\sigma(t_i)}$ within inter-switching time intervals which is block diagonalizable under the same similarity transformation and with the same block diagonal matrices sizes as the matrix $A_0(t_i)$, the stable (antistable) block diagonal matrix $\tilde{A}_0^-(t_i)(\tau-t_i)$ ($\tilde{A}_0^+(t_i)(\tau-t_i)$) generating a convergent (divergent) exponential matrix function $\tilde{\Phi}^-(\tau, t_i)$ ($\tilde{\Phi}^+(\tau, t_i)$).

**(iv)** Its first part follows from (2.8) since for any real constants $\alpha, \beta$ and any $A_{0j} \in \mathbf{A}_0$ which commute,

$$e^{A_{0j_1}\alpha} \cdot e^{A_{0j_2}\beta} = e^{\left(A_{0j_1}\alpha + A_{0j_2}\beta\right)}; \; \forall j_{1,2} \in \overline{N}_t.$$ Its second part follows from the semigroup property of $\Phi(t,\tau)$.

**(v)** It follows from (2.8) and (2.13), both being block diagonal with two non-zero square block matrices of corresponding identical time-invariant sizes, respectively $n^-$ and $n^+$, under the given assumptions since the matrices $\mathbf{A}_i \; \forall i \in \overline{q}$ are diagonalizable with identical two square matrix blocks of identical sizes.

If all the matrices in the set $\mathbf{A}_0$ are dichotomic, namely, they have no critically stable eigenvalues, then they admit a similarity transformation to a block diagonal form with only stable and instable eigenvalues. Under some extra assumptions related to the switching function to require a minimum residence time at each parameterization of the switched system, it may be proved that the input-state/ output operators of the solution are bounded operators. Now, denote by $\mathbf{P}_\pm^r : \mathbf{L}_2^r \to \mathbf{L}_{2\pm}^r$ the usual orthogonal projections of $\mathbf{L}_2^r$ onto $\mathbf{L}_{2\pm}^r$. Those projections are useful to describe the input-state and input-output operators for positive or negative times when the input is least-square integrable either for the negative or positive real semi-axis. The subsequent previous results is direct.

**Lemma 3.6**. $\Gamma_+ : \mathbf{L}_{2+}^m \to \mathbf{L}_{2+}^n$ and $\Gamma_{0+} : \mathbf{L}_{2+}^m \to \mathbf{L}_{2+}^p$ are linear bounded, and equivalently continuous, operators if any of the properties Theorem 2.1[(i)-(iii)] holds, and $\Gamma_+ : \mathbf{L}_{2-}^m \to \mathbf{L}_{2-}^n$ and $\Gamma_{0+} : \mathbf{L}_{2-}^m \to \mathbf{L}_{2-}^p$ are linear bounded, and equivalently continuous, operators if any of the properties in Theorem 2.2 holds.

**Proof**: It turns out from applying the Cauchy-Schwartz inequality to the sate/output - trajectory solutions that if the system is globally asymptotically stable and the input is an original (i.e. it identically zero for $t \in \mathbf{R}_-$) and, furthermore, square integrable, then the state and output trajectory solutions are identically zero for $t \in \mathbf{R}_-$ and square integrable on $\mathbf{R}_{0+}$. As a result, both linear operators are bounded and, equivalently, continuous. The second result is a dual one to the first result. □

Note that, compared to $\Gamma \in \mathbf{L}\left(\mathbf{L}_2^m, \mathbf{L}_2^n\right)$ and $\Gamma_o \in \mathbf{L}\left(\mathbf{L}_2^m, \mathbf{L}_2^p\right)$, the input-state operators $\Gamma_+ : \mathbf{L}_{2+}^m \to \mathbf{L}_{2+}^n$ and input-output $\Gamma_{o+} : \mathbf{L}_{2+}^m \to \mathbf{L}_{2+}^p$ (identified with the so-called causal Toeplitz operator if the input is an original vector function) have domains restricted from $\mathbf{L}_2^m$ to $\mathbf{L}_{2+}^m$ and projected images from $\mathbf{L}_2^n$, respectively $\mathbf{L}_2^p$, onto $\mathbf{L}_{2+}^n$, respectively $\mathbf{L}_{2+}^p$, provided that $\Gamma \in \mathbf{L}\left(\mathbf{L}_2^m, \mathbf{L}_2^n\right)$ and $\Gamma_o \in \mathbf{L}\left(\mathbf{L}_2^m, \mathbf{L}_2^p\right)$. In the same way, the input-state operators $\Gamma_- : \mathbf{L}_{2-}^m \to \mathbf{L}_{2-}^n$ and input-output $\Gamma_{o-} : \mathbf{L}_{2-}^m \to \mathbf{L}_{2-}^p$ have domains restricted from $\mathbf{L}_2^m$ to $\mathbf{L}_{2-}^m$ and projected images from $\mathbf{L}_2^n$, respectively $\mathbf{L}_2^p$, onto $\mathbf{L}_{2+}^n$, respectively $\mathbf{L}_{2+}^p$, provided that $\Gamma \in \mathbf{L}\left(\mathbf{L}_2^m, \mathbf{L}_2^n\right)$ and $\Gamma_o \in \mathbf{L}\left(\mathbf{L}_2^m, \mathbf{L}_2^p\right)$. Note also that Lemma 3.6 only gives sufficiency– type conditions of boundedness of those operators based on results of Theorems 2.1- 2.2. The following definitions are related to four important input to sate and input-output operators which are obtained from the operators $\Gamma \in \mathbf{L}\left(\mathbf{L}_2^m, \mathbf{L}_2^n\right)$ and $\Gamma_o \in \mathbf{L}\left(\mathbf{L}_2^m, \mathbf{L}_2^p\right)$ and which include domain restrictions and orthogonal projections of their images since they act on half axis Lebesgue spaces $\mathbf{L}_{2\pm}^{m,np}$:

**Definitions 3.7**. Let $\Gamma \in \mathbf{L}\left(\mathbf{L}_2^m, \mathbf{L}_2^n\right)$ bounded, so that $\Gamma_o \in \mathbf{L}\left(\mathbf{L}_2^m, \mathbf{L}_2^p\right)$, is also bounded. We define:

**1**. The causal input-output Hankel operator (or, simply causal Hankel operator) $\mathbf{H}_{\Gamma_o}: \mathbf{P}_+^p \Gamma_o |_{\mathbf{L}_{2-}^m} = \mathbf{P}_+^p \Gamma_o \mathbf{P}_-^m$ with symbol $\Gamma_o$.

**2**. The anticausal input-output Hankel operator (or, simply anticausal Hankel operator) $\hat{\mathbf{H}}_{\Gamma_o}: \mathbf{P}_-^p \Gamma_o |_{\mathbf{L}_{2+}^m} = \mathbf{P}_-^p \Gamma_o \mathbf{P}_+^m$ with symbol $\Gamma_o$.

**3**. The causal input-output Toeplitz operator (or, simply causal Toeplitz operator) $\mathbf{T}_{\Gamma_o}: \mathbf{P}_+^p \Gamma_o |_{\mathbf{L}_{2+}^m} = \mathbf{P}_+^p \Gamma_o \mathbf{P}_+^m$ with symbol $\Gamma_o$.

**4**. The anticausal input-output Toeplitz operator (or, simply anticausal Toeplitz operator) $\hat{\mathbf{T}}_{\Gamma_o}: \mathbf{P}_-^p \Gamma_o |_{\mathbf{L}_{2-}^m} = \mathbf{P}_-^p \Gamma_o \mathbf{P}_-^m$ with symbol $\Gamma_o$.

**5**. The causal input-state Hankel operator $\mathbf{H}_\Gamma: \mathbf{P}_+^n \Gamma |_{\mathbf{L}_{2-}^m} = \mathbf{P}_+^n \Gamma \mathbf{P}_-^m$ with symbol $\Gamma$.

**6**. The anticausal input-state Hankel operator $\hat{\mathbf{H}}_\Gamma: \mathbf{P}_-^n \Gamma |_{\mathbf{L}_{2+}^m} = \mathbf{P}_-^n \Gamma \mathbf{P}_+^m$ with symbol $\Gamma$.

**7**. The causal input-state Toeplitz operator $\mathbf{T}_\Gamma: \mathbf{P}_+^n \Gamma |_{\mathbf{L}_{2+}^m} = \mathbf{P}_+^n \Gamma \mathbf{P}_+^m$ with symbol $\Gamma$.

**8**. The anticausal input-state Toeplitz operator $\hat{\mathbf{T}}_\Gamma: \mathbf{P}_-^n \Gamma |_{\mathbf{L}_{2-}^m} = \mathbf{P}_-^n \Gamma \mathbf{P}_-^m$ with symbol $\Gamma$. □

The input - output Hankel and Toeplitz operators (see Definitions 3.7 [1-4] ), or simply Hankel and Toeplitz operators, are of wide use for the particular case of delay– free systems with single parameterizations, then being delay –free linear time- invariant systems (see, for instance, [29]). Definitions 3.7 and Lemmas 3.3-3.4 define extensions of those operators to describe the input-state/ output trajectories of the time-delayed switched system (2.1). The subsequent result related to the state and output trajectory solutions of the switched system (2.1) described by the input-sate and input-output Hankel and Toeplitz operators:

**Theorem 3.8**. The following properties hold under Assumption 3.1:

**(i)** $\mathbf{T}_{\Gamma_o} + \hat{\mathbf{H}}_{\Gamma_o} = \Gamma_o \mathbf{P}_+^m$, so that $\mathbf{T}_{\Gamma_o} = \mathbf{P}_+^p \Gamma_o \mathbf{P}_+^m$ iff $\hat{\mathbf{H}}_{\Gamma_0} = \mathbf{P}_-^p \Gamma_o \mathbf{P}_+^m = 0$, with

$$\left(\Gamma_o \mathbf{P}_+^m u\right)(t) = \int_0^t C(t)Z(t,\tau)B(\tau)u(\tau)d\tau + D(t)u(t)$$

$\hat{\mathbf{T}}_{\Gamma_o} + \mathbf{H}_{\Gamma_o} = \Gamma_o \mathbf{P}_-^m$, so that $\hat{\mathbf{T}}_{\Gamma_o} = \mathbf{P}_-^p \Gamma_o \mathbf{P}_-^m$ iff $\mathbf{H}_{\Gamma_o} = \mathbf{P}_+^p \Gamma_o \mathbf{P}_-^m = 0$, with

$$\left(\Gamma_o \mathbf{P}_-^m u\right)(t) = \int_{-\infty}^0 C(t)Z(t,\tau)B(\tau)u(\tau)d\tau + D(t)u(t)$$

**(ii)** $\mathbf{T}_\Gamma + \hat{\mathbf{H}}_\Gamma = \Gamma \mathbf{P}_+^m$, so that $\mathbf{T}_\Gamma = \mathbf{P}_+^p \Gamma \mathbf{P}_+^m$ iff $\hat{\mathbf{H}}_\Gamma = \mathbf{P}_-^p \Gamma \mathbf{P}_+^m = 0$, with

$$\left(\Gamma \mathbf{P}_+^m u\right)(t) = \int_0^t Z(t,\tau)B(\tau)u(\tau)d\tau$$

$\hat{T}_\Gamma + H_\Gamma = \Gamma P_-^m$, so that $\hat{T}_\Gamma = = P_-^p \Gamma P_-^m$ iff $H_\Gamma = P_+^p \Gamma P_-^m = 0$, with

$$\left(\Gamma P_-^m u\right)(t) = \int_{-\infty}^0 Z(t,\tau) B(\tau) u(\tau) d\tau$$

**(iii)** $\left(\hat{H}_{\Gamma_o} u\right)(t) = 0$

$$\left(T_{\Gamma_o} u\right)(t) = \left(P_+^p \Gamma_o P_-^m u\right)(t) = \int_0^t C(t) Z(t,\tau) B(\tau) u(\tau) d\tau + D(t) u(t) = \int_0^t \tilde{C}(t) \tilde{Z}(t,\tau) \tilde{B}(\tau) u(\tau) d\tau + \tilde{D}(t) u(t)$$

$$= \sum_{j=1}^{k(t)} \int_{t_{j-1}}^{t_j} C(t_{j-1}) \left[ \Phi(t,\tau) B(t_{j-1}) + \sum_{i=1}^q \sum_{\ell=1}^{k(t)} \int_{t_{\ell-1}}^{t_\ell} \Phi(t,\gamma) A_i(t_{\ell-1}) Z(\gamma - h_i, \tau) B(t_{j-1}) (U(\tau) - U(\gamma - h_i)) d\gamma \right] u(\tau) d\tau$$

$$+ D(t) u(t)$$

$$= \sum_{j=1}^{k(t)} \int_{t_{j-1}}^{t_j} \tilde{C}(t_{j-1}) \left[ \tilde{\Phi}(t,\tau) \tilde{B}(t_{j-1}) + \sum_{i=1}^q \sum_{\ell=1}^{k(t)} \int_{t_{\ell-1}}^{t_\ell} \tilde{\Phi}(t,\gamma) \tilde{A}_i(t_{\ell-1}) \tilde{Z}(\gamma - h_i, \tau) \tilde{B}(t_{j-1}) (U(\tau) - U(\gamma - h_i)) d\gamma \right] u(\tau) d\tau$$

$$+ \tilde{D}(t) u(t) \quad ; \forall t \in \mathbf{R}_{0+}$$

The last expression being valid if $t_1 = 0$ since $t_0 = -\infty$. If $t_1 > 0$ then the given switching sequence $ST_\sigma$ may be redefined as $t_1 \to 0, t_{i+1} \to t_i, \ldots \forall i \geq 1$ with the switching function initialized as $\sigma(t) = \sigma(t_0) = -\infty; \forall t \in (-\infty, t_2]$ so that the switched system is not modified and the above expression is valid for the causal Toeplitz operator.

$$\left(T_{\Gamma_o} u\right)(t) = 0, \quad \left(\hat{H}_{\Gamma_o} u\right)(t) = \left(P_-^p \Gamma_o P_+^m u\right)(t) = -\int_0^\infty \tilde{C}^+(t) \left(\tilde{Z}(t,\tau) \tilde{B}(\tau)\right)^+ u(\tau) d\tau$$

$$= -\sum_{j=k(t)}^\infty \int_{t_j}^{t_{j+1}} \tilde{C}^+(t) \left[ \tilde{\Phi}^+(t,\tau) \tilde{B}^+(t_j) + \sum_{i=1}^q \sum_{\ell=1}^{k(t)} \int_{t_{\ell-1}}^{t_\ell} \left(\tilde{\Phi}(t,\gamma) \tilde{A}_i(t_{\ell-1}) \tilde{Z}(\gamma - h_i, \tau) \tilde{B}(t_j)\right)^+ (U(\tau) - U(\gamma - h_i)) d\gamma \right] u(\tau) d\tau$$

$$; \forall t \in \mathbf{R}_{0-}$$

with the switching time instants being redefined with $t_1 = 0$, so that $\sigma(t) = \sigma(-\infty); \forall t \in (-\infty, t_2]$, as above, in the case that the first switching time instant is nonzero.

**(iv)** $\left(\hat{T}_{\Gamma_o} u\right)(t) = 0$

$$\left(H_{\Gamma_o} u\right)(t) = \left(P_+^p \Gamma_o P_-^m u\right)(t) = \int_{-\infty}^0 C(t) Z(t,\tau) B(\tau) u(\tau) d\tau = \int_{-\infty}^0 \tilde{C}(t) \tilde{Z}(t,\tau) \tilde{B}(\tau) u(\tau) d\tau$$

$$= \int_{-\infty}^0 C(-\infty) \left[ \Phi(t,\tau) B(-\infty) + \sum_{i=1}^q \int_{-\infty}^0 \Phi(t,\gamma) A_i(-\infty) Z(\gamma - h_i, \tau) B(-\infty) (U(\tau) - U(\gamma - h_i)) d\gamma \right] u(\tau) d\tau$$

$$= \int_{-\infty}^0 \tilde{C}(-\infty) \left[ \tilde{\Phi}(t,\tau) \tilde{B}(-\infty) + \sum_{i=1}^q \int_{-\infty}^0 \tilde{\Phi}(t,\gamma) \tilde{A}_i(-\infty) \tilde{Z}(\gamma - h_i, \tau) \tilde{B}(-\infty) (U(\tau) - U(\gamma - h_i)) d\gamma \right] u(\tau) d\tau$$

$$; \forall t \in \mathbf{R}_{0+}$$

$\left(H_{\Gamma_o} u\right)(t) = 0$

$$\left(\hat{\mathbf{T}}_{\Gamma_o}u\right)(t)=\left(\mathbf{P}_-^p\Gamma\mathbf{P}_-^m u\right)(t)=\int_{-\infty}^0 C(t)Z(t,\tau)B(\tau)u(\tau)d\tau+D(t)u(t)=\int_{-\infty}^0 \tilde{C}(t)\tilde{Z}(t,\tau)\tilde{B}(\tau)u(\tau)d\tau+\tilde{D}(t)u(t)$$

$$=\int_{-\infty}^0 C^-(-\infty)\left[\Phi^-(t,\tau)B^-(-\infty)+\sum_{i=1}^q\int_{-\infty}^0\left(\Phi(t,\gamma)A_i(-\infty)Z(\gamma-h_i,\tau)B(-\infty)\right)^-(U(\tau)-U(\gamma-h_i))d\gamma\right]u(\tau)d\tau$$

$$+D(t)u(t)$$

$$=\int_{-\infty}^0 \tilde{C}^-(-\infty)\left[\tilde{\Phi}^-(t,\tau)\tilde{B}^-(-\infty)+\sum_{i=1}^q\int_{-\infty}^0\left(\tilde{\Phi}(t,\gamma)\tilde{A}_i(-\infty)\tilde{Z}(\gamma-h_i,\tau)\tilde{B}(-\infty)\right)^-(U(\tau)-U(\gamma-h_i))d\gamma\right]u(\tau)d\tau$$

$$+\tilde{D}(t)u(t) \qquad\qquad ;\ \forall t\in\mathbf{R}_{0-}$$

(v) $\left(\hat{\mathbf{H}}_\Gamma u\right)(t)=0$

$$(\mathbf{T}_\Gamma u)(t)=\left(\mathbf{P}_+^n\Gamma\mathbf{P}_-^m u\right)(t)=\int_0^t Z(t,\tau)B(\tau)u(\tau)d\tau=\int_0^t \tilde{Z}(t,\tau)\tilde{B}(\tau)u(\tau)d\tau$$

$$=\sum_{j=1}^{k(t)}\int_{t_{j-1}}^{t_j}\left[\Phi(t,\tau)B(t_{j-1})+\sum_{i=1}^q\sum_{\ell=1}^{k(t)}\int_{t_{\ell-1}}^{t_\ell}\Phi(t,\gamma)A_i(t_{\ell-1})Z(\gamma-h_i,\tau)B(t_{j-1})(U(\tau)-U(\gamma-h_i))d\gamma\right]u(\tau)d\tau$$

$$=\sum_{j=1}^{k(t)}\left[\tilde{\Phi}(t,\tau)\tilde{B}(t_{j-1})+\sum_{i=1}^q\sum_{\ell=1}^{k(t)}\int_{t_{\ell-1}}^{t_\ell}\tilde{\Phi}(t,\gamma)\tilde{A}_i(t_{\ell-1})\tilde{Z}(\gamma-h_i,\tau)\tilde{B}(t_{j-1})(U(\tau)-U(\gamma-h_i))d\gamma\right]u(\tau)d\tau$$

$$;\ \forall t\in\mathbf{R}_{0+}$$

The last expression being valid if $t_1=0$ since $t_0=-\infty$. If $t_1>0$ then the given switching sequence $ST_\sigma$ may be redefined as $t_1\to 0, t_{i+1}\to t_i,\ldots \forall i\geq 1$ with the switching function initialized as $\sigma(t)=\sigma(t_0)=-\infty;\forall t\in(-\infty,t_2]$ so that the switched system is not modified and the above expression is+ valid for the causal input-state Toeplitz operator.

$$\left(\mathbf{T}_\Gamma u\right)(t)=0\ ,\ \left(\hat{\mathbf{H}}_\Gamma u\right)(t)=\left(\mathbf{P}_-^n\Gamma\mathbf{P}_+^m u\right)(t)=-\int_0^\infty\left(\tilde{Z}(t,\tau)\tilde{B}(\tau)\right)^+u(\tau)d\tau$$

$$=-\sum_{j=k(t)}^\infty\int_{t_j}^{t_{j+1}}\left[\tilde{\Phi}^+(t,\tau)\tilde{B}^+(t_j)+\sum_{i=1}^q\sum_{\ell=1}^{k(t)}\int_{t_{\ell-1}}^{t_\ell}\left(\tilde{\Phi}(t,\gamma)\tilde{A}_i(t_{\ell-1})\tilde{Z}(\gamma-h_i,\tau)\tilde{B}(t_j)\right)^+(U(\tau)-U(\gamma-h_i))d\gamma\right]u(\tau)d\tau$$

$$;\ \forall t\in\mathbf{R}_{0-}$$

with the switching time instants being redefined with $t_1=0$, so that $\sigma(t)=\sigma(t_0)=\sigma(-\infty);\forall t\in(-\infty,t_2]$, as above, in the case that the first switching time instant is nonzero.

(vi) $\left(\hat{\mathbf{T}}_\Gamma u\right)(t)=0$

$$\left(\mathbf{H}_\Gamma u\right)(t)=\left(\mathbf{P}_+^n\Gamma_0\mathbf{P}_-^m u\right)(t)=\int_{-\infty}^0 Z(t,\tau)B(\tau)u(\tau)d=\int_{-\infty}^0 \tilde{Z}(t,\tau)\tilde{B}(\tau)u(\tau)d\tau$$

$$=\int_{-\infty}^0 C(-\infty)\left[\Phi(t,\tau)B(-\infty)+\sum_{i=1}^q\int_{-\infty}^0\Phi(t,\gamma)A_i(-\infty)Z(\gamma-h_i,\tau)B(-\infty)(U(\tau)-U(\gamma-h_i))d\gamma\right]u(\tau)d\tau$$

$$= \int_{-\infty}^{0} \tilde{C}(-\infty) \left[ \tilde{\Phi}(t,\tau)\tilde{B}(-\infty) + \sum_{i=1}^{q} \int_{-\infty}^{0} \tilde{\Phi}(t,\gamma)\tilde{A}_i(-\infty)\tilde{Z}(\gamma-h_i,\tau)\tilde{B}(-\infty)(U(\tau)-U(\gamma-h_i))d\gamma \right] u(\tau) d\tau$$

$$; \forall t \in \mathbf{R}_{0+}$$

$$(\mathbf{H}_\Gamma u)(t) = 0$$

$$(\hat{\mathbf{T}}_\Gamma u)(t) = (\mathbf{P}_-^n \Gamma \mathbf{P}_-^m u)(t) = \int_{-\infty}^{0} Z(t,\tau)B(\tau)u(\tau)d\tau = \int_{-\infty}^{0} \tilde{Z}(t,\tau)\tilde{B}(\tau)u(\tau)d\tau$$

$$= \int_{-\infty}^{0} C(-\infty)\left[\Phi(t,\tau)B(-\infty) + \sum_{i=1}^{q} \int_{-\infty}^{0} \Phi(t,\gamma)A_i(-\infty)Z(\gamma-h_i,\tau)B(-\infty)(U(\tau)-U(\gamma-h_i))d\gamma\right] u(\tau)d\tau$$

$$= \int_{-\infty}^{0} \tilde{C}(-\infty)\left[\tilde{\Phi}(t,\tau)\tilde{B}^-(-\infty) + \sum_{i=1}^{q}\int_{-\infty}^{0}\tilde{\Phi}(t,\gamma)\tilde{A}_i(-\infty)\tilde{Z}(\gamma-h_i,\tau)\tilde{B}(-\infty)(U(\tau)-U(\gamma-h_i))d\gamma\right]u(\tau)d\tau$$

$$; \forall t \in \mathbf{R}_{0-}$$

**Proof**: It follows directly from Lemmas 3.3-3.4 and Definitions 3.7 by noting that $\Gamma \in \mathbf{L}\left(\mathbf{L}_2^m, \mathbf{L}_2^n\right)$ and $\Gamma_o \in \mathbf{L}\left(\mathbf{L}_2^m, \mathbf{L}_2^p\right)$ are bounded operators from Assumption 3.1 since all configurations of the switched system have no critically stable eigenvalues and, furthermore, no stable or unstable ones within the open vertical strip $(-\varepsilon, \varepsilon) \times \mathbf{R} \subset \mathbf{C}$ from Theorems 2.1-2.2. □

Note that, if Assumptions 3.2 hold, then Theorem (3.8) holds with a constant transformation of coordinates $T \in \mathbf{R}^{n \times n}$ in (3.4)-(3.5), i.e.

$$\tilde{A}(t) = \tilde{A}_{\sigma(t_k)} = TA(t)T^{-1} = TA_{\sigma(t_k)}T^{-1}$$

$$\tilde{B}(t) = \tilde{B}_{\sigma(t_k)} = TB(t) = TB_{\sigma(t_k)} \quad, \quad \tilde{C}(t) = \tilde{C}_{\sigma(t_k)} = C(t)T^{-1} = C_{\sigma(t_k)}T^{-1}$$

$\forall t \in [t_k, t_{k+1}), \forall t_k, t_{k+1} \in ST_\sigma, \forall t \in [t_\ell, \infty)$ if $t_\ell \in ST_\sigma$ and there is no $ST_\sigma \ni t > t_\ell$ so that $\text{card}(ST_\sigma) < \infty$. Theorem 3.8 can be specified as follows under Assumptions 3.8 provided that the matrices of delayed dynamics have each two block diagonal expressions of the same orders as those of $\mathbf{A}_0$.

**Corollary 3.9**. Assume that all the matrices in the set of configurations $\mathbf{A}_i$ ($\forall i \in \overline{q} \cup \{0\}$) are block diagonal with two matrix blocks matrices of orders $n_-$ and $n_+$ identical to those of the stable and antistable blocks of the matrices in the set $\mathbf{A}_0$ consisting of the N of delay – free matrices of dynamics. Thus, Theorem 3.8 has the following particular expressions for the anticausal (input-output) Hankel and input-state Hankel operators provided that Assumptions 3.2 hold:

$$(\hat{\mathbf{H}}_\Gamma u)(t) = -\int_0^\infty \tilde{Z}^+(t,\tau)\tilde{B}(\tau)^+ u(\tau)d\tau$$

$$= -\sum_{j=k(t)}^{\infty} \int_{t_j}^{t_{j+1}} \left[ \tilde{\Phi}^+(t,\tau)\tilde{B}^+(t_j) + \sum_{i=1}^{q}\sum_{\ell=1}^{k(t)} \int_{t_{\ell-1}}^{t_\ell} \tilde{\Phi}^+(t,\gamma)\tilde{A}_i^+(t_{\ell-1})\tilde{Z}^+(\gamma-h_i,\tau)\tilde{B}^+(t_j)(U(\tau)-U(\gamma-h_i))d\gamma \right] u(\tau)d\tau$$

$$; \forall t \in \mathbf{R}_{0-}$$

so that $(\hat{H}_{\Gamma_o} u)(t) = \tilde{C}^+(t)(\hat{H}_\Gamma u)(t); \forall t \in R_{0-}$, with the switching time instants being redefined with $t_1 = 0$, so that $\sigma(t) = \sigma(t_0) = \sigma(-\infty); \forall t \in (-\infty, t_2]$, as above, in the case that the first switching time instant is nonzero.

**Proof**: It follows directly from Theorem 3.8 and Lemma 3.5 (iii) –(iv) since the matrix functions $\Phi(t,\tau)$ and $Z(t,\tau)$ maintain a two block diagonal structure with matrices of orders $n_-$ and $n_+$ from (3.10). □

Definitions of causality and anticausality follow:

**Definition 3.10**. A bounded input-output linear operator $\Gamma_o : L_2^m \to L_2^p$ is said to be causal (anticausal) if the anticausal Hankel operator is zero, i.e. $\hat{H}_{\Gamma_o} = 0$ (if the causal Hankel operator is zero, i.e. $H_{\Gamma_o} = 0$). □

**Definition 3.11**. The switched system (2.1) is said to be causal (anticausal) if $\hat{H}_{\Gamma_o} = 0$ ($H_{\Gamma_o} = 0$) provided that $\Gamma_o : L_2^m \to L_2^p$ is bounded. □

**Definition 3.12**. A bounded input-state linear operator $\Gamma : L_2^m \to L_2^n$ is said to be causal (anticausal) if $\hat{H}_\Gamma = 0$ ($H_\Gamma = 0$). □

A direct result from Definitions 3.10-3.12 is the following :

**Assertion 3.13**. If $\Gamma : L_2^m \to L_2^n$ is bounded and causal (anticausal) then the switched system (2.1) is causal (anticausal) but the converse is not true, in general. □

**Theorem 3.14**. The following properties hold under Assumptions 3.1 for a given switching function $\sigma : R_{0+} \to \overline{N}$ provided that it obeys a minimum residence time between consecutive switches which exceeds some appropriate minimum threshold:

**(i)** $\Gamma_o : L_2^m \to L_2^p$ is bounded independent of the delays and if all the matrices of delay-free dynamics in the set $A_0$ are stable then the system (2.1) is globally asymptotically stable and causal independent of the delays.

**(ii)** If $\hat{H}_{\Gamma_o} : L_2^m \to L_2^p$ is zero independent of the delays and the switched system (2.1) is uniformly controllable and uniformly observable independent of the delays then it is globally asymptotically Lyapunov´s stable independent of the delays.

**(iii)** If $\hat{H}_\Gamma : L_2^m \to L_2^n$ is zero independent of the delays and the switched system (2.1) is uniformly controllable independent of the delays then it is globally asymptotically Lyapunov´s stable independent of the delays.

**Proof**: (i) $\Gamma_o : L_2^m \to L_2^p$ is bounded from Assumptions 3.1, Theorem 2.1 and Theorem 2.2 if there is a sufficiently large residence time for the given switching function since there is an eigenvalue-free open vertical strip including the imaginary complex axis for any delays. Thus, all the configurations of the switched system are dichotomic independent of the delays if the switching function is subject to a minimum residence time exceeding an appropriate threshold. From Theorem 3.8 (i) and (iv), the system is causal if the anticausal Hankel operator is zero, namely,

$$\left(\hat{\mathbf{H}}_{\Gamma_o} u\right)(t) = -\int_0^\infty \tilde{C}^+(t)\left(\tilde{Z}(t,\tau)\tilde{B}(\tau)\right)^+ u(\tau)d\tau = 0 \;;\; \forall u \in L_{2+}^m, \forall t \in \mathbf{R}_{0-}$$

Property **(i)** follows since if $\mathbf{A}_0$ is a set of stable matrices then the switched system is globally asymptotically stable independent of the delays from Theorem 2.21 and Assumptions 3.1 and causal from $\left(\hat{\mathbf{H}}_{\Gamma_o} u\right)(t) = 0$; $\forall u \in L_{2+}^m, \forall t \in \mathbf{R}_{0-}$. The above factorization exists since $0 \le n_+(t) < \infty$ (the number of unstable eigenvalues of any configuration of (2.1) is finite) since the characteristic quasi-polynomials of all the configurations have a principal term in view of the structure of (2.1), [17]. Since the system is uniformly observable, then the following contradiction is stated if $n_+(t) \ne 0$:

$$\left(\hat{\mathbf{H}}_{\Gamma_o} u\right)(t) = \int_0^\infty \tilde{C}^+(t)\left(\tilde{Z}(t,\tau)\tilde{B}(\tau)\right)^+ u(\tau)d\tau = 0$$

$$\Rightarrow 0 = \int_0^\infty \left(\tilde{Z}(t,\tau)\tilde{B}(\tau)\right)^+ u(\tau)d\tau \ge \sum_{i=1}^\ell \int_{t_i}^{t_{i+1}} \left(\tilde{Z}(t,\tau)\tilde{B}(\tau)\right)^+ \left(\tilde{Z}(t,\tau)\tilde{B}(\tau)\right)^{+T} g(t_i)d\tau \ne 0, \forall t \in \mathbf{R}_{0-}$$

provided that the control $u(\tau) = \left(\tilde{Z}(t,\tau)\tilde{B}(\tau)\right)^{+T} g(t_i)$, $\forall \tau \in [t_i, t_{i+1})$, $\forall t_i, t_{i+1} \in ST_\sigma$ for $\left\{0 \ne \mathbf{R}^n \ni g(t_i) = o\left(1/\left\|\left(\tilde{Z}(t,\tau)\tilde{B}(\tau)\right)\right\|_{L_2[t_i, t_{i+1}]}\right)\right\}_0^\infty$. The contradiction follows since $t_{i+1} \ge t_i + T$ for some $T \in \mathbf{R}_{0+}$ so that the controllability Grammian $\int_{t_i}^{t_i + T_0} \left(\tilde{Z}(t,\tau)\tilde{B}(\tau)\right)\left(\tilde{Z}(t,\tau)\tilde{B}(\tau)\right)^{+T} d\tau$, $\forall \tau \in [t_i, t_{i+1})$ is positive definite, $\forall t_i, t_{i+1} \in ST_\sigma$ for some constant $T_0 \in \mathbf{R}_{0+}$ (independent of $t_i$) and $\forall t_i \in ST_\sigma$ if the system (2.1) is uniformly controllable. Thus, $n_+(t) = 0$ and Property (ii) follows. Property (iii) follows in a similar way by neglecting the controllability condition since

$$\left(\hat{\mathbf{H}}_\Gamma u\right)(t) = \int_0^\infty \left(\tilde{Z}(t,\tau)\tilde{B}(\tau)\right)^+ u(\tau)d\tau = 0 \,. \qquad \square$$

The following result strengths Theorem 3-14 since Assumptions 3.2 allow to maintain all the eigenvalues strictly outside the imaginary axis independent of the delays via arbitrary switching (see Theorem 2.1 and Theorem 2.2):

**Corollary 3.15**. If Assumptions 3.2 hold then Theorem 3.14 holds for an arbitrary switching function, i.e. without requiring a minimum residence time in-between any two consecutive active parameterizations.

$\square$

Theorem 3.14 has the following simpler version for zero and small delays which follows from the continuity of the eigenvalues with respect to the delays. It is not required that the matrices describing the delayed dynamics of the various configurations have sufficiently small norms compared with the minimum absolute stability abscissa among the configurations associated with the delay – free dynamics defined by the set $\mathbf{A}_0$:

**Theorem 3.16.** Assume that:

a) the set of Assumptions 3.1 holds except the stability conditions in Theorem 2.1 (i.e. there is no requirement on the smallness of the norms of the matrices describing the delayed dynamics of the various configurations of the switched system)

b) a switching function $\sigma: \mathbf{R}_{0+} \to \overline{N}$ is given which respects a minimum residence time between consecutive switches exceeding some appropriate minimum threshold

Then, the following properties hold:

(i) If $\mathbf{\Gamma}_o : \mathbf{L}_2^m \to \mathbf{L}_2^p$ is bounded for some switching function $\sigma: \mathbf{R}_{0+} \to \overline{N}$ and if all the matrices of delay-free dynamics in the set $\mathbf{A}_0$ are stable then the system (2.1) is globally asymptotically stable and causal for $h_i \in [0, \overline{h}), \forall i \in \overline{q}$ for some sufficiently small $\overline{h} \in \mathbf{R}_+$.

(ii) If $\hat{\mathbf{H}}_{\Gamma_o} : \mathbf{L}_2^m \to \mathbf{L}_2^p$ is zero and the switched system (2.1) is uniformly controllable and uniformly observable for $h_i \in [0, \overline{h}), \forall i \in \overline{q}$ for some sufficiently small $\overline{h} \in \mathbf{R}_+$ then it is globally asymptotically Lyapunov´s stable for $h_i \in [0, \overline{h}), \forall i \in \overline{q}$.

(iii) If $\hat{\mathbf{H}}_{\Gamma} : \mathbf{L}_2^m \to \mathbf{L}_2^n$ is zero and the switched (2.1) is uniformly controllable for $h_i \in [0, \overline{h}), \forall i \in \overline{q}$ for some sufficiently small $\overline{h} \in \mathbf{R}_+$ then it is globally asymptotically Lyapunov´s stable for $h_i \in [0, \overline{h}), \forall i \in \overline{q}$. □

The following result follows from Theorem 3.16 under Assumptions 3.2 in the same way as Corollary 3.15 is a consequence of Theorem 3.14:

**Corollary 3.17.** If Assumptions 3.2 hold then Theorem 3.16 holds for an arbitrary switching function, i.e. without requiring a minimum residence time in-between any two consecutive active parameterizations.

□

The condition of the auxiliary unforced delay – free system being dichotomic can be removed to conclude global asymptotic stability under causality and uniform controllability and observability as proved in the sequel:

**Corollary 3.18.** If the switched system (2.1) is causal, uniformly controllable and uniformly observable independent of the delays for a given switching function then it is globally asymptotically stable independent of the delays.

**Proof**: Define the truncated linear operator $\mathbf{\Gamma}_{ot} : \mathbf{L}_{2t}^m \to \mathbf{L}_{2t}^p$ for arbitrary (but finite) $t \in \mathbf{R}$ in the same way as $\mathbf{\Gamma}_o : \mathbf{L}_2^m \to \mathbf{L}_2^p$ where $\mathbf{L}_{2t}^r := \{ f \in \mathbf{L}_2^r : f(\tau) = 0, \forall \tau \in (-\infty, t) \cup (t, \infty) \}$ is the set of square - integrable r- real vector functions of compact support $[-t, t]$. Note that the linear operator $\mathbf{\Gamma}_{ot}$ is bounded for any finite time " t" irrespective of the spectrum of the system. If the system is causal, uniformly observable and uniformly controllable then for some point non-singular n- transformation matrix function $T(t,.):[t, t+T_0] \to \mathbf{R}^{n \times n}$ [see Lemma 3.5 (i)], a matrix function $\tilde{Z}(t,\tau)$ being similar to $Z(t,\tau)$ may be calculated leading to:

$$(\hat{H}_{\Gamma_o} u)(t) = 0 \Rightarrow$$

$$(\hat{H}_\Gamma u)(t) = \int_0^\infty (\tilde{Z}(t,\tau)\tilde{B}(\tau))^+ (\tilde{Z}(t,\tau)\tilde{B}(\tau))^{+T} u(\tau) d\tau = 0 \geq \int_t^{t+T_0} (\tilde{Z}(t,\tau)\tilde{B}(\tau))^+ (\tilde{Z}(t,\tau)\tilde{B}(\tau))^{+T} g(t) d\tau = 0$$

for some any finite $t \in \mathbf{R}_{0+}$, some constant $T_0 \in \mathbf{R}_+$ and some $g \in \mathbf{L}_{2,t+T_0} \supset \mathbf{L}_{2,t+\delta}$; $\forall \delta \in [0, T_0]$, chosen so that $0 \neq u(\tau) = (\tilde{Z}(t,\tau)\tilde{B}(\tau))^{+T} g \in \mathbf{L}_{2,t+T_0}$; $\forall \tau \in [t, t+T_0]$. The superscript "+" now includes the contribution of the finite number of unstable and critically unstable modes (since the system is not assumed to be dichotomic) and the integrand is a square matrix function of piecewise constant order $n^+ : [t, t+T_0] \to \mathbf{Z}_{0+}$; for any finite real t. Such a matrix order function is finite, since the whole number of critically stable and unstable modes is always finite since all the configurations of the switched system have a principal term in its characteristic quasi-polynomial. This, together with the finiteness of $T_0$ implies that the controllability Grammian $\int_t^{t+T_0} (\tilde{Z}(t,\tau)\tilde{B}(\tau))^+ (\tilde{Z}(t,\tau)\tilde{B}(\tau))^{+T} g(t) d\tau$ may be decomposed in a finite sum of matrices of constant order $\bar{n}^+(t) := \max(n^+(\tau) : \tau \in [t, t+T_0])$ completed if necessary with zero blocks for the remaining terms in the sum, the number of additive terms being the number of discontinuities in $n^+(\tau)$ plus one. This leads again to a contradiction as in Theorem 3.14 and the causal system being uniformly controllable and observable cannot possess critically unstable and unstable modes. □

ACKNOWLEDGMENTS

The authors are very grateful to the Spanish Ministry of Education by its partial support of this work through project DPI2006-00714. They are also grateful to the Basque Government by its support through GIC07143-IT-269-07, SAIOTEK SPED06UN10 and SPE07UN04.